\magnification=\magstep1
\input amstex
\documentstyle{amsppt}

\font\tencyr=wncyr10 
\font\sevencyr=wncyr7 
\font\fivecyr=wncyr5 
\newfam\cyrfam \textfont\cyrfam=\tencyr \scriptfont\cyrfam=\sevencyr
\scriptscriptfont\cyrfam=\fivecyr
\define\hexfam#1{\ifcase\number#1
  0\or 1\or 2\or 3\or 4\or 5\or 6\or 7 \or
  8\or 9\or A\or B\or C\or D\or E\or F \or G\or H\fi}
\mathchardef\Sha="0\hexfam\cyrfam 58

\define\defeq{\overset{\text{def}}\to=}

\def \isom {\buildrel \sim \over \rightarrow}
\def \mosi {\buildrel \sim \over \leftarrow}

\def \Im{\operatorname {Im}}
\def \pr{\operatorname {pr}}
\def \tor{\operatorname {tor}}
\def \cl{\operatorname {cl}}
\def \id{\operatorname {id}}

\def \Ker{\operatorname {Ker}}
\def \char{\operatorname {char}}
\def \Spec{\operatorname {Spec}}

\def \Gal{\operatorname {Gal}}
\def \cl{\operatorname {cl}}
\def \Hom{\operatorname {Hom}}
\def \res{\operatorname {res}}
\def \inf{\operatorname {inf}}
\def \Sel{\operatorname {Sel}}
\def \Coker{\operatorname {Coker}}

\def \Tr{\operatorname {Tr}}
\def \Aut{\operatorname {Aut}}
\def \id{\operatorname {id}}
\def \Ext{\operatorname {Ext}}
\def \GL{\operatorname {GL}}
\def \Mor{\operatorname {Mor}}
\def \ab{\operatorname {ab}}

\def \prof{\operatorname {prof}}
\def \prol{\text{pro-$l$}}
\def \sep{\operatorname{sep}}
\def \wild{\operatorname{w}}
\def \tame{\operatorname{t}}
\def \et{\text{\'et}}
\def \cpt{\operatorname {cpt}}
\def \diag{\operatorname {diag}}
\def \Coker{\operatorname {Coker}}
\define\nemp{\neq\emptyset}

\define\Primes{\frak{Primes}}
\define\Se{\frak{Sel}}
\define\Sh{\frak{Sha}}

\define\ZSigma{\widehat{\Bbb Z}^{\Sigma}}

\NoRunningHeads
\NoBlackBoxes
\topmatter

\title
On the arithmetic of abelian varieties
\endtitle
\bigskip

\author
{MOHAMED SA\"IDI and AKIO TAMAGAWA}
\endauthor

\abstract 
We prove some new results on the arithmetic of abelian varieties over function fields of 
one variable over finitely generated (infinite) fields. Among other things, we introduce certain 
new natural objects `discrete Selmer groups' and `discrete Shafarevich-Tate groups', 
and prove that they are finitely generated $\Bbb Z$-modules. Further, we prove that 
in the isotrivial case, the discrete Shafarevich-Tate group vanishes and the discrete 
Selmer group coincides with the Mordell-Weil group. 
One of the key ingredients to prove these results is a new
specialisation theorem 
for first Galois cohomology groups,  
which generalises N\'eron's specialisation theorem
for rational points of abelian varieties. 
\endabstract

\toc

\subhead
\S0. Introduction
\endsubhead

\subhead
\S1. A Specialisation Theorem for $H^1$
\endsubhead

\subhead
\S2. Selmer Groups
\endsubhead

\subhead
\S3. Shafarevich-Tate Groups
\endsubhead

\subhead
\S4. The Isotrivial Case
\endsubhead

\endtoc

\endtopmatter
\document

\subhead
\S 0. Introduction
\endsubhead
Let $k$ be a field of characteristic $0$, and $C\to \Spec k$ 
a smooth, separated and geometrically connected (not necessarily proper)
{\it algebraic curve} over $k$. Write $K=k(C)$ for the function field of $C$, 
$C^{\cl}$ for the set of closed points of $C$, 
and $k(c)$ for the residue field of $C$ at 
$c\in C^{\cl}$. 
Let 
$\Cal A\to C$ 
be an {\it abelian scheme} over $C$. Write $A\defeq \Cal A_K \defeq \Cal A \times _C \Spec K$ 
for the {\it generic fibre} 
of $\Cal A$. 
For 
each $c\in C^{\cl}$, write $\Cal A_c\defeq \Cal A\times _C \Spec k(c)$ for the fibre of $\Cal A$ at $c$,  
$K_c$ for the completion of $K$ at $c$, and $A_c\defeq A\times _KK_c$. 
Thus, $A$ (resp. $\Cal A_c$, resp. $A_c$) is an abelian variety over $K$ (resp. $k(c)$, 
resp. $K_c$). 
Consider the Kummer exact sequence
$$0 \to A(K)^{\wedge} \to H^1(G_K,TA) \to TH^1(G_K,A) \to 0,$$ 
where $TA$ is the (full) Tate module of $A$, $A(K)^{\wedge}$ is the completion 
$\varprojlim_{N>0} A(K)/NA(K)$ of the group $A(K)$ of 
$K$-rational points of $A$ (which coincides with the profinite completion of $A(K)$, 
if $A(K)$ is finitely generated), and  $TH^1(G_K,A)$ is the (full) Tate module of the Galois 
cohomology group $H^1(G_K,A)$ classifying $K$-principal homogeneous spaces under $A$.
Similarly, for each closed point $c\in C^{\cl}$, we have the Kummer exact sequences
$$0 \to A_c(K_c)^{\wedge} \to H^1(G_{K_c},T A_c) \to TH^1(G_{K_c},A_c) \to 0$$ 
and 
$$0 \to \Cal A_c(k(c))^{\wedge} \to H^1(G_{k(c)},T\Cal A_c) \to TH^1(G_{k(c)},\Cal A_c) \to 0.$$ 

We have a natural commutative diagram

$$
\CD
0 @>>> A(K)^{\wedge} @>>> H^1(G_K,TA) @>>> TH^1(G_K,A) @>>> 0  \\
@. @VVV      @V
VV    @V
VV \\
0 @>>> \prod_cA_c(K_c)^{\wedge} @>>>  \prod_c H^1(G_{K_c},T A_c) @>>> \prod_c TH^1(G_{K_c},A_c) @>>> 0 \\
\endCD
\tag {0.1}
$$

\noindent
where the horizontal sequences are the above Kummer exact sequences, the vertical maps are natural 
restriction maps, and the product is taken over all closed points $c\in C^{\cl}$. In fact, 
when $k$ is finitely generated over $\Bbb Q$, diagram (0.1) can be identified with the following 
natural commutative diagram (cf. Proposition 2.1 (ii)): 

$$
\CD
0 @>>> \Cal A(C)^{\wedge} @>>> H^1(\pi_1(C),TA) @>>> TH^1_{\et}(C,\Cal A) @>>> 0  \\
@. @VVV      @V
VV    @V
VV \\
0 @>>> \prod_c \Cal A_c(k(c))^{\wedge} @>>>  \prod_c H^1(G_{k(c)},T \Cal A_c) @>>> 
\prod_c TH^1(G_{k(c)},\Cal A_c) @>>> 0 \\
\endCD
\tag {0.2}
$$

\noindent
where the upper horizontal sequence is a Kummer exact sequence for the \'etale site of $C$, 
the lower horizontal sequence is as above, 
the vertical maps are natural 
restriction maps, 
and the product is taken over all closed points $c\in C^{\cl}$. 

Just as in the case where $K$ is a number field, define the {\it profinite Selmer group} 
$$\Sel (A)\defeq \Sel (A,C)\defeq \Ker \lgroup H^1(G_K,TA) \to \prod_c TH^1(G_{K_c},A_c)\rgroup,$$
and the {\it Shafarevich-Tate group}
$$\Sha(A)\defeq \Sha (A,C)\defeq \Ker \lgroup H^1(G_K,A)\to \prod_c H^1(G_{K_c},A_c)\rgroup.$$
Thus, we have a natural exact sequence
$$0\to A(K)^{\wedge}\to \Sel(A)\to T\Sha(A)\to 0,$$
where $T\Sha(A)$ is the Tate module of $\Sha(A)$.
For an integer $N>0$, define the {\it $N$-Selmer group} by
$$\Sel_N(A)\defeq \Sel (A,C)_N\defeq \Ker \lgroup H^1(G_K,A[N])\to \prod _c H^1(G_{K_c},A_c)\rgroup,$$ 
so that $\Sel (A)=\varprojlim_{N>0}\Sel_N(A)$. 

One of our main results is the following, which improves a result of [Lang-Tate] 
(cf. Proposition 2.10, Proposition 3.9 (i) and Remark 3.13). 

\proclaim {Proposition A} Assume that $k$ is {\bf finitely generated} over $\Bbb Q$. 
Then for each integer $N>0$, the $N$-Selmer group
$\Sel_N(A)$, as well as the subgroup $\Sha(A)[N]$ of $N$-torsion points of $\Sha(A)$, is {\bf finite}. 
\endproclaim

The proof of Proposition A follows from the following specialisation result (cf. Proposition 1.8).

\proclaim {Proposition B} Assume that $k$ is {\bf Hilbertian} (cf [Serre2], 9.5).
Then for each integer $N>0$, 
there exists a {\bf finite} subset $S\subset C^{\cl}$ (of cardinality $\leq 2$), 
depending on $N$, such that the natural restriction map 
$$H^1(\pi_1(C),A[N]) @>>> \prod_{c\in S} H^1(G_{k(c)},\Cal A_c[N]),$$ 
is {\bf injective}. 
\endproclaim

We also prove the following analogous 
specialisation result for the Galois cohomology of the $l$-adic Tate module of $A$.
(cf. Proposition 1.4).
\proclaim {Proposition C} Assume that $k$ is {\bf Hilbertian}. Let $l$ be a prime number. 
Then there exists a {\bf finite} subset 
$S\subset C^{\cl}$ of cardinality $\leq 2$, 
depending on $l$, such that the natural restriction map
$$H^1(\pi_1(C),T_lA)\to \prod_{c\in S} H^1(G_{k(c)},T_l\Cal A_c)$$ 
is {\bf injective}. 
\endproclaim

In the case where either $k$ is {\bf finitely generated} over $\Bbb Q$ 
or the $\overline k$-trace of $A_{K\overline k}\defeq A\times_KK\overline k$ is trivial, 
one can prove that there exists a finite subset $S\subset C^{\cl}$ as in Proposition C
{\bf of cardinality $1$} (cf. Proposition 1.2 and Proposition 1.4). 
We do not know 
(even in the finitely generated case) 
if an analogue of Proposition C holds for the Galois cohomology of the full Tate module $TA$.

As a consequence of Proposition C, one deduces the following (cf. Proposition 2.2).

\proclaim{Proposition D} Assume that $k$ is {\bf Hilbertian}. 
Then the middle and left vertical maps in diagrams (0.1) and (0.2) are {\bf injective}.
\endproclaim

For the rest of this introduction we will assume that $k$ is {\bf finitely generated} over $\Bbb Q$.
We will identify $A(K)^{\wedge}$, $H^1(G_K,TA)$, and $\prod_c\Cal A_c(k(c))^{\wedge}$ with their images in
$\prod_c H^1(G_{k(c)},T\Cal A_c)$. For each closed point $c\in C^{\cl}$ the group $\Cal A_c(k(c))$ of $k(c)$-rational points of $\Cal A_c$
is finitely generated as $k(c)$ is finitely generated over $\Bbb Q$ (Mordell-Weil Theorem, cf. [Lang-N\'eron]), hence injects into 
its profinite completion
$\Cal A_c(k(c))^{\wedge}$. We identify $\Cal A_c(k(c))$ with its image in $\Cal A_c(k(c))^{\wedge}$.
We define the {\bf discrete} Selmer group by
$$\Se (A)\defeq \Se (A,C)\defeq \Sel (A) \bigcap \prod_c \Cal A_c(k(c)) \subset \prod_cH^1(G_{k(c)},T\Cal A_c).$$
Note that $A(K)\subset \Se(A)$. 
We define the {\bf discrete} Shafarevich-Tate group by 
$$\Sh(A)\defeq\Sh(A,C)\defeq\Se(A)/A(K).$$

We conjecture the following (cf. Conjecture 3.8).

\definition {Conjecture E}
The equality $\Se(A)=A(K)$ (or, equivalently, $\Sh(A)=0$) holds. 
\enddefinition

Concerning Conjecture E, we prove the following (cf. Proposition 2.5, Proposition 3.3 and Proposition 3.7).

\proclaim {Proposition F} The discrete Selmer group $\Se(A)$ is a {\bf finitely generated} $\Bbb Z$-module. 
The discrete Shafarevich-Tate group $\Sh(A)$ is a 
{\bf finitely generated free} $\Bbb Z$-module. 
\endproclaim

\proclaim {Proposition G}
Assume that there exists a prime number $l$ such that the $l$-primary part $\Sha(A)[l^{\infty}]$ of the torsion group
$\Sha(A)$ is finite. 
Then the assertion of 
Conjecture E holds. 
\endproclaim

Our results are most complete in the case where the abelian variety $A$ is {\it isotrivial}. In this case we prove the following (cf. Theorem 4.1).

\proclaim {Theorem H} Assume that 
the abelian variety $A$ is {\bf isotrivial}, i.e., 
$A_{\overline K}$ descends to an abelian variety 
over $\overline k$. 
Then the Shafarevich-Tate group $\Sha (A)$ is {\bf finite}. 
In particular, the assertion of Conjecture E holds in this case. 
\endproclaim

Although we assumed above that $\char(k)=0$, we prove similar results in arbitrary characteristics. 

Some of the results in this paper have applications in anabelian geometry. More precisely, 
Conjecture E and Theorem H have applications to Grothendieck's anabelian section conjecture 
(cf. [Sa\"\i di1], $\S0$, for a precise statement of this conjecture). 
One can prove that the validity of Conjecture E above implies that the section conjecture 
(for $\pi_1$ of proper hyperbolic curves) over finitely generated fields reduces to the case 
of number fields. Using among others Theorem H, one can also prove that if the section conjecture holds 
for all proper hyperbolic curves over all number fields then it holds for all proper hyperbolic curves 
over all finitely generated fields which are defined over a number field (cf. [Sa\"\i di2], $\S5$).

Finally, we explain the content of each section briefly. 
In $\S1$, we prove Propositions B and C.
In $\S2$ and $\S3$, we prove Propositions A, D, F and G. 
In $\S4$, we prove Theorem H. 

\subhead
Notations
\endsubhead
Next, we fix notations that will be used throughout this paper.

Given a (profinite) group $G$ and a (continuous) $G$-module $C$, we write $C^G\defeq H^0(G,C)$.

Let $H$ be an abelian group. 
For an integer $N>0$, we write 
$H/N\defeq H/NH$ and $H[N]\defeq \{h\in H\ \vert\ Nh=0\}$. 
We write 
$H^{\wedge}\defeq
\varprojlim_{N>0} H/N$, and 
$H^{\prof}\defeq 
\varprojlim_{H'\subset H,\  (H:H')<\infty}
H/H'$ 
for the profinite completion of $H$. 
Thus, we have natural homomorphisms 
$H\otimes _{\Bbb Z}\widehat {\Bbb Z}\to 
H^{\wedge}\to H^{\prof}$, which are isomorphisms when $H$ is finitely generated. 
We write $H^{\tor}\defeq \bigcup _{N>0}H[N]$ 
for the torsion subgroup of $H$, and set $H/\{\tor\}\defeq H/H^{\tor}$.

For a prime number $l$, 
we write 
$H^{\wedge,l}\defeq \varprojlim_{n\geq 0} H/l^n$ for the 
$l$-adic completion of $H$, and 
$H^{\prol}\defeq 
\varprojlim_{H'\subset H,\  (H:H'): \text{$l$-power}}
H/H'$ 
for the pro-$l$ completion of $H$. 
Thus, we have natural homomorphisms 
$H\otimes _{\Bbb Z}\Bbb Z_l\to 
H^{\wedge}\otimes _{\widehat{\Bbb Z}}\Bbb Z_l
=H^{\wedge,l}\to H^{\prol}$, which are isomorphisms when $H$ is finitely generated. 
We write $H^{\tor,l}$ and $H^{\tor,l'}$ for the $l$-primary part and the prime-to-$l$ part, 
respectively, of the torsion abelian group $H^{\tor}$, and set 
$H/\{\text{$l$-tor}\}\defeq H/H^{\tor,l}$ and 
$H/\{\text{$l'$-tor}\}\defeq H/H^{\tor,l'}$. 

Let $\Primes$ be the set of all prime numbers. For a nonempty subset $\Sigma \subset \Primes$, 
we say that $N$ is a $\Sigma$-integer if $N$ is divisible only by primes in $\Sigma$. 
We set 
$\ZSigma 
\defeq 
\varprojlim_{N:\text{$\Sigma$-integer$>0$}} 
\Bbb Z/N =
\prod_{l\in \Sigma}\Bbb Z_l$, 
write 
$H^{\wedge,\Sigma}\defeq \varprojlim_{N:\text{$\Sigma$-integer$>0$}} 
H/N
= \prod_{l\in \Sigma} H^{\wedge,l}$ 
(thus, $H^{\wedge,\{l\}}= H^{\wedge,l}$) 
for the $\Sigma$-adic completion of $H$, 
and write 
$T^{\Sigma}H\defeq \varprojlim_{N:\text{$\Sigma$-integer$>0$}} 
H[N]=
\prod_{l\in \Sigma}T_l H$ (where 
$T_l H\defeq T^{\{l\}}H$) 
for the $\Sigma$-adic Tate module of $H$. 
(Note that $T^{\Sigma}H$ is always torsion-free.) 
We write $H^{\tor,\Sigma}\defeq \bigcup _{\text{$N$:$\Sigma$-integer$>0$}}H[N]
=\oplus_{l\in\Sigma}H^{\tor,l}$ 
for the $\Sigma$-primary torsion subgroup of $H$, and set $H/\{\text{$\Sigma$-tor}\}\defeq H/H^{\tor,\Sigma}$. 
Set $\Sigma'\defeq \Primes\setminus\Sigma$. 

Let $B$ be an abelian variety 
over a field $\kappa$ of characteristic $p\ge 0$ 
with algebraic closure $\overline \kappa$ and 
separable closure $\kappa^{\sep}\subset\overline\kappa$. We write 
$B[N]$, 
$B^{\tor}$, 
$T^{\Sigma}B$, 
$T_lB$ 
instead of 
$B(\kappa^{\sep})[N]$, 
$B(\kappa^{\sep})^{\tor}$, 
$T^{\Sigma}B(\kappa^{\sep})$, 
$T_lB(\kappa^{\sep})$, 
respectively. 
Write $\Primes^{\dag}\defeq \Primes\setminus \{p\}$. 
For a nonempty subset 
$\Sigma \subset \Primes^{\dag}$, 
recall the Kummer exact sequence
$$0\to B(\kappa)^{\wedge,\Sigma}\to H^1(G_{\kappa},T^{\Sigma}B)\to T^{\Sigma}H^1(G_{\kappa},B)\to 0,\tag 0.3$$
where $G_{\kappa}\defeq \Gal (\kappa^{\sep}/\kappa)$ 
is the absolute Galois group of $\kappa$ and  
$H^1(G_{\kappa},B)\defeq H^1(G_{\kappa}, B(\kappa^{\sep}))$ 
is the first Galois cohomology group, 
which arises from the Kummer exact sequence of $G_{\kappa}$-modules 
$$0\to B
[N]\to B(\kappa^{\sep}) @>N>> B(\kappa^{\sep})\to 0,$$ 
where $N$ denotes the map of multiplication by a $\Sigma$-integer $N>0$. 
Note that  the above sequence (0.3) induces a natural isomorphism
$(B(\kappa)^{\wedge,\Sigma})^{\tor}\isom H^1(G_{\kappa},T^{\Sigma}B)^{\tor}$, 
as $T^{\Sigma}H^1(G_{\kappa},B)$ is torsion-free. 


\subhead
\S1. A Specialisation Theorem for $H^1$
\endsubhead
Let $S$ be a locally noetherian, regular, integral scheme. 
Write $K$ for the function field of $S$, 
$k(t)$ for the residue field of $S$ at $t\in S$, and 
$p_t  (\geq 0)$ for the characteristic of $k(t)$.  
Write $\char(S)\defeq\{p_t\mid t\in S\} \subset \Primes \cup \{0\}$. 
Let 
$\Cal A\to S$
be an {\it abelian scheme} over $S$. We write $A\defeq \Cal A_K \defeq \Cal A \times _S \Spec K$ for the 
{\it generic fibre} 
of $\Cal A$, and, for each $t\in S$, 
we write $\Cal A_t\defeq \Cal A\times _S \Spec k(t)$ for the fibre of $\Cal A$ at $t$. 
Thus, $A$ (resp. $\Cal A_t$) is an abelian variety over $K$ (resp. over $k(t)$).

Let $\eta$ be a geometric point of $S$ with values in the generic point of $S$. 
Then $\eta$ determines an algebraic 
closure 
$\overline K$ 
and a separable 
closure $K^{\sep}$ 
of $K$, 
Write $G_K=\Gal(K^{\sep}/K)$ 
for the absolute Galois group of 
$K$, and 
$\pi_1(S)=\pi_1(S,\eta)$ for the 
\'etale fundamental group of $S$. 
Thus, we have a natural exact sequence of profinite groups
$$1\to I_S\to G_K\to \pi_1(S)\to 1$$
where $I_S$ is defined so that the sequence is exact. 

Write $S^1$ for the set of points of codimension $1$ of $S$.  
For each $t\in S^1$, the local ring 
$\Cal O_{S,t}$ is a discrete valuation ring, and 
let $(G_K\supset) D_t\supset I_t$ be a {\it decomposition group} and an {\it inertia group} 
associated to $t$. Thus, $D_t$ and $I_t$ are only defined up to conjugation in $G_K$. 
We have a natural exact sequence
$$ 1\to I_t\to D_t\to G_{k(t)}\to 1$$
where $G_{k(t)}\defeq \Gal (k(t)^{\sep}/k(t))$. 
Then, by purity, the group $I_S$ 
is (topologically) normally generated by the subgroups $I_t$, where $t$ runs over all points in $S^1$. 
We have a natural exact sequence 
$$1\to I_t^{\wild}\to I_t\to I_t^{\tame}\to 1,$$
where the wild inertia group $I_t^{\wild}$ is defined to be the unique Sylow-$p_t$ subgroup of $I_t$ 
(resp. the trivial subgroup $\{1\}\subset I_t$) for $p_t>0$ (resp. $p_t=0$), and the tame inertia 
group $I_t^{\tame}$ is defined by $I_t^{\tame}\defeq I_t/I_t^{\wild}$. 
Note that $I_t^{\tame}$ is naturally isomorphic to $\widehat {\Bbb Z}^{p_t'}(1)$,
where 
$\widehat {\Bbb Z}^{p_t'}\defeq\widehat{\Bbb Z}
^{\Primes \setminus\{p_t\}}
$, 
and the ``$(1)$'' denotes a Tate twist.

\proclaim {Lemma 1.1} Let $\Sigma \subset \Primes \setminus\char(S)$ be a nonempty subset. Then: 

\noindent
{\rm (i)} The $G_K$-module $T^{\Sigma}A$ (hence, in particular, 
$\prod_{l\in \Sigma}A[l]$) has a natural structure of $\pi_1(S)$-module. 

\noindent
{\rm (ii)} For each $l\in\Sigma$, 
write 
$\pi_1(S)[A,l]$ for the kernel of the natural map 
$\pi_1(S)\to \Aut(A[l])$ (cf. (i)), 
$\pi_1(S)[A,l]^l$ for the maximal pro-$l$ quotient of 
$\pi_1(S)[A,l]$, and 
$\pi_1(S)(A,l)$ for the kernel of the natural surjection 
$\pi_1(S)[A,l]\twoheadrightarrow \pi_1(S)[A,l]^l$. 
Define 
$\pi_1(S)[A,\Sigma]\defeq \bigcap _{l\in \Sigma} \pi_1(S)[A,l]$,  
and $\pi_1(S)(A,\Sigma)\defeq \bigcap _{l\in \Sigma}\pi_1(S)(A,l)$, where the intersection is over all prime integers $l\in \Sigma$. 
Further, let $\Pi_S^{A,\Sigma}\defeq \pi_1(S)/ \pi_1(S)(A,\Sigma)$. (Note that  $\pi_1(S)(A,\Sigma)$ is a normal subgroup of 
$\pi_1(S)$ since $\pi_1(S)(A,l)$ is a characteristic subgroup of $\pi_1(S)[A,l]$.)
Thus, we have a natural exact sequence 
$$1\to \pi_1(S)[A,\Sigma]/\pi_1(S)(A,\Sigma) \to \Pi_S^{A,\Sigma} \to \pi_1(S)/\pi_1(S)[A,\Sigma]\to 1,$$
where 
$$\pi_1(S)/\pi_1(S)[A,\Sigma]
\hookrightarrow \prod_{l\in\Sigma}(\pi_1(S)/\pi_1(S)[A,l])
\hookrightarrow \prod_{l\in\Sigma}\Aut(A[l])$$ 
and 
$$\pi_1(S)[A,\Sigma]/\pi_1(S)(A,\Sigma)=\prod_{l\in\Sigma}\Im(\pi_1(S)[A,\Sigma]\to \pi_1(S)[A,l]^l).$$ 
Then the $\pi_1(S)$-module $T^{\Sigma}A$ has a natural structure of 
$\Pi_S^{A,\Sigma}$-module. 

\noindent
{\rm (iii)} 
The natural inflation map 
$H^1(\Pi_S^{A,\Sigma},T^{\Sigma}A) \to H^1(\pi_1(S),T^{\Sigma}A)$ 
is an isomorphism. The natural inflation map 
$H^1(\pi_1(S),T^{\Sigma}A) \to H^1(G_K,T^{\Sigma}A)$ 
is an injection in general, and an isomorphism 
if $k(t)$ is finitely generated over the prime field 
for each $t\in S^1$. 

\noindent
{\rm (iv)} Let $N$ be a $\Sigma$-integer $>0$. Then 
the natural inflation map 
$H^1(\Pi_S^{A,\Sigma},A[N]) \to H^1(\pi_1(S),A[N])$ 
is an isomorphism and the natural inflation map 
$H^1(\pi_1(S),A[N]) \to H^1(G_K,A[N])$ 
is an injection. 
\endproclaim

\demo {Proof} 
(i) For each $t\in S^1$, any inertia group $I_t$ associated to $t$ acts trivially
on $T^{\Sigma}A$, as follows from the well-known criterion of 
good reduction for abelian varieties (cf. [Serre-Tate], Theorem 1). 
Thus, $I_S$ acts trivially on $T^{\Sigma}A$, 
and $T^{\Sigma}A$ has a natural structure of $\pi_1(S)$-module. 

\noindent
(ii) This follows from (i) and the fact that, for each $l\in\Sigma$,  
$\Ker(\Aut(T_lA) \to \Aut(A[l]))$ 
($\simeq 
\Ker(\GL_{2d}(\Bbb Z_l) \to \GL_{2d}(\Bbb F_l))$, 
where $d\defeq\dim A$) is pro-$l$. 

\noindent
(iii) First, note that the inflation maps for various $H^1$ are always injective, 
and that various $H^1$ with coefficients in $T^{\Sigma}A$ decompose into the 
direct product of $H^1$ with coefficients in $T_lA$ for $l\in\Sigma$. Thus, 
it suffices to prove, for each $l\in\Sigma$, that the inflation map 
$H^1(\Pi_S^{A,l}, T_lA)\to H^1(\pi_1(S),T_lA)$ 
is an isomorphism  in general 
(where $\Pi_S^{A,l}\defeq \Pi_S^{A,\{l\}}$), 
and that 
the inflation map 
$H^1(\pi_1(S),T_lA)\to H^1(G_K, T_lA)$ is an isomorphism 
if $k(t)$ is finitely generated over the prime field 
for each $t\in S^1$. 

For the first assertion, 
consider the 
inflation-restriction exact sequence
$$0\to H^1(\Pi_S^{A,l},T_lA) @>\inf>> H^1(\pi_1(S),T_lA) @>\res>> 
H^1(\pi_1(S)(A,l),T_lA)^{\Pi_S^{A,l}}.$$
We claim that $H^1 (\pi_1(S)(A,l),T_lA)=0$. 
Indeed, this follows from 
the fact that 
$H^1 (\pi_1(S)(A,l),T_lA)=\Hom(\pi_1(S)(A,l),T_lA)$, and 
$\pi_1(S)(A,l)=\Ker(\pi_1(S)[A,l]\twoheadrightarrow \pi_1(S)[A,l]^l)$. 

For the second assertion, 
consider the 
inflation-restriction exact sequence
$$0\to H^1(\pi_1(S),T_lA) @>\inf>> H^1(G_K,T_lA) @>\res>> 
H^1(I_S,T_lA)^{\pi_1(S)}.$$
We claim that $H^1(I_S,T_lA)^{\pi_1(S)}=0$. 
Indeed, this follows from the fact that 
$H^1 (I_S,T_lA)=\Hom(I_S,T_lA)$ 
and that, for each $t\in S^1$, 
one has 
$\Hom (I_t,T_lA)^{G_{k(t)}}=\Hom(I_t^{\tame,l}, T_lA)^{G_{k(t)}}=0$, 
where $I_t^{\tame,l}$ is the maximal pro-$l$ quotient of $I_t^{\tame}$. 
This last statement follows from the fact that the action of $G_{k(t)}$ 
on $I_t^{\tame,l}\isom \Bbb Z_l(1)$ and $T_lA$ has different weights. 
More precisely, $T_l\Cal A_t \simeq 
T_lA$ as $G_{k(t)}$-modules, and 
the $G_{k(t)}$-representation
$I_t^{\tame,l}\otimes \Bbb Q_l$ (resp. $T_l{\Cal A_t}\otimes \Bbb Q_l$) is pure of weight $-2$ (resp. pure of weight $-1$) since $k(t)$ is finitey generated (cf. [Jannsen], 2). 
Thus, $\Hom(I_t^{\tame,l}, T_l{\Cal A_t})^{G_{k(t)}}=0$ follows (cf. loc. cit., Fact 2).

\noindent
(iv) As the inflation maps for various $H^1$ are always injective 
and various $H^1$ with coefficients in $A[N]$ decompose into the 
direct product of $H^1$ with coefficients in $A[l^r]$ for $l\in\Sigma$ and 
some $r\geq 0$, it suffices to prove, for each $l\in\Sigma$, 
that the inflation map 
$H^1(\Pi_S^{A,l}, A[l^r])\to H^1(\pi_1(S),A[l^r])$ 
is an isomorphism. For this, 
consider the 
inflation-restriction exact sequence
$$0\to H^1(\Pi_S^{A,l},A[l^r]) @>\inf>> H^1(\pi_1(S),A[l^r]) @>\res>> 
H^1(\pi_1(S)(A,l),A[l^r])^{\Pi_S^{A,l}}.$$
We claim that $H^1 (\pi_1(S)(A,l),A[l^r])=0$. 
Indeed, this follows from 
the fact that 
$H^1 (\pi_1(S)(A,l),A[l^r])=\Hom(\pi_1(S)(A,l),A[l^r])$, 
and 
$\pi_1(S)(A,l)=\Ker(\pi_1(S)[A,l]\twoheadrightarrow \pi_1(S)[A,l]^l)$. 
\qed
\enddemo

Next, 
let $k$ be a field of characteristic $p\ge 0$, and $C\to \Spec k$ a smooth, separated 
and geometrically connected (not necessarily proper)
{\it algebraic curve} over $k$. 
Write $K=k(C)$ for the function field of $C$, 
$C^{\cl}$ for the set of closed points of $C$ 
(which coincides with the set $C^1$ of codimension $1$ of $C$), 
and $k(c)$ for the residue field of $C$ at $c\in C^{\cl}$.  Let 
$\Cal A\to C$
be an {\it abelian scheme} over $C$. We write $A\defeq \Cal A_K \defeq \Cal A \times _C \Spec K$ for the 
{\it generic fibre} 
of $\Cal A$, and, for each $c\in C^{\cl}$, 
we write $\Cal A_c\defeq \Cal A\times _C \Spec k(c)$ for the fibre of $\Cal A$ at $c$. 
Thus, $A$ (resp. $\Cal A_c$) is an abelian variety over $K$ (resp. over $k(c)$).

Let $\eta$ be a geometric point of $C$ with values in the generic point of $C$. 
Then $\eta$ determines algebraic closures $\overline K$ and $\overline k$ 
and separable closures $K^{\sep}$ and $k^{\sep}$ 
of $K$ and $k$, respectively, and a geometric point $\overline \eta$ of 
$ C_{\overline k}\defeq C\times _k\overline k$. 
Write $G_k=\Gal(k^{\sep}/k)$, $G_K=\Gal(K^{\sep}/K)$ and 
$G_{Kk^{\sep}}=\Gal(K^{\sep}/Kk^{\sep})$ for the absolute Galois groups of 
$k$, $K$ and $Kk^{\sep}$, respectively. Write 
$\pi_1(C)=\pi_1(C,\eta)$ and $\pi_1( C_{\overline k})=\pi_1( C_{\overline k},\overline\eta)$ for the 
\'etale fundamental groups of $C$ and $ C_{\overline k}$, respectively. 
Thus, we have natural exact sequences of profinite groups
$$1\to \pi_1( C_{\overline k})\to \pi_1(C)\to G_k \to 1,\tag 1.1$$
$$1\to G_{Kk^{\sep}} \to G_K\to G_k\to 1,\tag 1.2$$
and
$$1\to I_C\to G_K\to \pi_1(C)\to 1, \tag 1.3$$
where $I_C$ is defined so that sequence (1.3) is exact. 
We have a commutative diagram of exact sequences

$$
\CD
@. 1 @. 1\\
@. @VVV   @VVV\\
@.  I_C  @=  I_C\\
@. @VVV  @VVV\\
1@>>> G_{Kk^{\sep}} @>>> G_K@>>>G_k@>>> 1 \\
@. @VVV   @VVV   \Vert\\
1@>>> \pi_1( C_{\overline k})@>>> \pi_1(C)@>>> G_k @>>> 1\\
@. @VVV @VVV\\
@. 1 @. 1\\
\endCD
$$ 

For each $c\in C^{\cl}$, write $K_c$ (resp. $\widehat{\Cal O}_{C,c}$) for the completion of $K$ 
(resp. $\Cal O_{C,c}$) at $c$, and 
$A_c\defeq A\times_KK_c$. 
Thus, $K_c$ (resp. $\widehat{\Cal O}_{C,c}$) is a complete discrete
valuation field (resp. ring) of equal characteristic $p\ge 0$ with residue field $k(c)$, and 
$A_c$ is an abelian variety over $K_c$. 
Let $(G_K\supset)\  D_c\supset I_c$ be a {\it decomposition group} and an {\it inertia group} 
associated to $c$. Thus, $D_c$ and $I_c$ are only defined up to conjugation in $G_K$. We have a natural exact sequence
$$ 1\to I_c\to D_c\to G_{k(c)}\to 1, \tag 1.4$$
where $G_{k(c)}\defeq \Gal (k(c)^{\sep}/k(c))$ is identified with the image of $D_c$ in $G_k$ (cf. sequence (1.2)). 
Then the group $I_C$ (cf. sequence (1.3))
is (topologically) normally generated by the subgroups $I_c$, where $c$ runs over all points in $C^{\cl}$. 
We have a natural exact sequence 
$$1\to I_c^{\wild}\to I_c\to I_c^{\tame}\to 1,$$
where the wild inertia group $I_c^{\wild}$ is defined to be the unique Sylow-$p$ subgroup of $I_c$ 
(resp. the trivial subgroup $\{1\}\subset I_c$) for $p>0$ (resp. $p=0$), and the tame inertia 
group $I_c^{\tame}$ is defined by $I_c^{\tame}\defeq I_c/I_c^{\wild}$. 
Note that $I_c^{\tame}$ is naturally isomorphic to $\widehat {\Bbb Z}^{\dag}(1)$,
where $\Primes ^{\dag}\defeq \Primes \setminus\{p\}$, 
$\widehat {\Bbb Z}^{\dag}\defeq\widehat{\Bbb Z}^{\Primes^{\dag}}$, 
and the ``$(1)$'' denotes a Tate twist. 


For each $c\in C^{\cl}$, we have a natural commutative diagram (up to conjugation): 
$$
\matrix
G_{K_c}&\isom&D_c&\subset &G_K\\
&&&&\\
&&\downarrow&&\downarrow\\
&&&&\\
&&G_{k(c)}&\overset{s_c}\to{\hookrightarrow}& \pi_1(C)
\endmatrix
$$ 
where the vertical arrows are natural surjections and the horizontal arrows are natural injections. 
(The map $s_c: G_{k(c)}=\pi_1(\Spec k(c)) \to \pi_1(C)$ is associated to the natural morphism 
$\Spec k(c)\to C$ with image $c$ by functoriality of $\pi_1$.) Further, this diagram induces 
natural commutative diagrams 
$$
\matrix
H^1(\pi_1(C), A[N]) &\hookrightarrow& H^1(G_K, A[N])\\
&&\\
\downarrow&&\downarrow\\
&&\\
H^1(G_{k(c)}, \Cal A_c[N]) &\hookrightarrow& H^1(G_{K_c}, A_c[N])\\
\endmatrix
$$ 
for each $\Primes^{\dag}$-integer $N>0$, and 
$$
\matrix
H^1(\pi_1(C), T^{\Sigma}A) &\hookrightarrow& H^1(G_K, T^{\Sigma}A)\\
&&\\
\downarrow&&\downarrow\\
&&\\
H^1(G_{k(c)}, T^{\Sigma}\Cal A_c) &\hookrightarrow& H^1(G_{K_c}, T^{\Sigma}A_c)\\
\endmatrix
$$ 
for each nonempty subset $\Sigma\subset\Primes^{\dag}$, where 
the horizontal arrows are inflation maps and the vertical arrows are 
natural restriction maps. 


One of our main results in this section is the following.

\proclaim{Proposition 1.2} Let $\Sigma\subset \Primes^{\dag}$ be a 
{\bf finite} subset. 
Assume that $k$ is {\bf finitely generated} over the prime field and infinite.
Then there exists a closed point $c\defeq c(\Sigma)\in C^{\cl}$, 
depending on ($A$ and) $\Sigma$, such that the natural restriction map
$
H^1(\pi_1(C),T^{\Sigma}A) \to H^1(G_{k(c)},T^{\Sigma}\Cal A_c)$ is {\bf injective}.
\endproclaim

\demo{Proof} 
For simplicity, write $Q$ for the prime field of $k$ and $Z$ for the image of $\Bbb Z[\frac{1}{l}\ ; \ l\in\Sigma]$ in $Q$. 
Thus, we have $Q=\Bbb Q$ (resp. $Q=\Bbb F_p$) and $Z=\Bbb Z[\frac{1}{l}\ ; \ l\in\Sigma]$ (resp. $Z=\Bbb F_p$) when 
$p=0$ (resp. $p>0$). Then, as $k$ is finitely generated over the perfect field $Q$, the system 
$\Cal A \to C \to \Spec k\to\Spec Q$ admits a smooth model 
$\widetilde{\Cal A}\to\Cal C \to V \to U$. More precisely, $U=\Spec Z$; $V$ is an integral scheme which is smooth over $U$ 
and whose function field is isomorphic to (and is identified with) $k$; $\Cal C$ is a smooth scheme over $V$ whose 
generic fibre $\Cal C\times_V k$ is $k$-isomorphic to (and is identified with) $C$; and 
$\widetilde{\Cal A}\to\Cal C$ is an abelian scheme such that $\widetilde{\Cal A}\times_{\Cal C}C$ is isomorphic to 
(and is identified with) $\Cal A\to C$. 
Let $\pi_1(\Cal C)$ be the fundamental group of $\Cal C$, with respect to the base point which is induced by the base point $\eta$
in the beginning of $\S1$. Thus, there exists a natural continuous surjective homomorphism 
$\pi_1(C)\twoheadrightarrow \pi_1(\Cal C)\twoheadrightarrow \Pi_{\Cal C}^{A,\Sigma}$. 
By Lemma 1.1 (ii), $T^\Sigma A$ has a natural structure of 
$\Pi_{\Cal C}^{A,\Sigma}$-module, and, by Lemma 1.1 (iii), 
the inflation maps 
$$H^1(\Pi_{\Cal C}^{A,\Sigma}, T^\Sigma A)
\to H^1(\pi_1(\Cal C),T^\Sigma A) 
\to H^1(\pi_1(C),T^\Sigma A)\to H^1(G_K, T^\Sigma A)$$
are all isomorphisms. Recall that 
we have a natural exact sequence 
$$1\to \pi_1(\Cal C)[A,\Sigma]/\pi_1(\Cal C)(A,\Sigma) \to \Pi_{\Cal C}^{A,\Sigma} \to \pi_1(\Cal C)/\pi_1(\Cal C)[A,\Sigma]\to 1,$$
where 
$$\pi_1(\Cal C)/\pi_1(\Cal C)[A,\Sigma]
\hookrightarrow \prod_{l\in\Sigma}(\pi_1(\Cal C)/\pi_1(\Cal C)[A,l])
\hookrightarrow \prod_{l\in\Sigma}\Aut(A[l])$$ 
and 
$$\pi_1(\Cal C)[A,\Sigma]/\pi_1(\Cal C)(A,\Sigma)=\prod_{l\in\Sigma}\Im(\pi_1(\Cal C)[A,\Sigma]\to \pi_1(\Cal C)[A,l]^l).$$ 
In particular, as $\Sigma$ is finite, $\pi_1(\Cal C)/\pi_1(\Cal C)[A,\Sigma]$ is finite. 

\bigskip\noindent
{\bf Claim:} The Frattini subgroup of $\Pi_{\Cal C}^{A,\Sigma}$ is open. 

\bigskip
Indeed, as $\pi_1(\Cal C)/\pi_1(\Cal C)[A,\Sigma]$ is finite 
and $\pi_1(\Cal C)[A,\Sigma]/\pi_1(\Cal C)(A,\Sigma)$ is a direct product of pro-$l$ groups for $l\in\Sigma$, 
it suffices to prove 
that 
$\pi_1(\Cal C)[A,\Sigma]^l$ is a finitely generated pro-$l$ group for each $l\in\Sigma$ (cf. [Serre2], 10.6, Proposition), 
or, equivalently, 
that $\pi_1(\Cal C)[A,\Sigma]^{\ab}/l$ is finite. Let $\Cal C'\to \Cal C$ be the finite \'etale cover 
corresponding to the open subgroup 
$\pi_1(\Cal C)[A,\Sigma]\subset \pi_1(\Cal C)$, so that 
$\pi_1(\Cal C)[A,\Sigma]$ is identified with $\pi_1(\Cal C')$. 
Now, the desired finiteness follows from [Katz-Lang]. More precisely, 
let $Q'$ be the algebraic closure of $Q$ in $k$, which is finite over $Q$, and 
$Z'$ the integral closure of $Z$ in $Q'$. Then the morphism 
$\Cal C'\to\Spec Z$ 
factors as 
$\Cal C'\to\Spec Z'\to\Spec Z$, and the morphism 
$\Cal C'\times_ZQ \to\Spec Q$ factors as 
$\Cal C'\times_ZQ \to\Spec Q'\to\Spec Q$. 
The morphism $\Cal C'\times_ZQ \to\Spec Q'$ is smooth and geometrically connected, 
hence there exist an open subscheme 
$\Cal C'_1\subset\Cal C'$ (containing $\Cal C'\times_QQ'$) and an open subscheme 
$W\subset \Spec Z'$, such that 
$\Cal C'\to\Spec Z'$ induces a smooth, surjective morphism $\Cal C'_1\to W$ 
of finite type with geometrically connected generic fibre. 
Now, by [Katz-Lang], Lemma 2 (2), we have an exact sequence 
$$0\to \Ker(\Cal C'_1/W)\to \pi_1(\Cal C'_1)^{\ab}\to \pi_1(W)^{\ab}\to 0,$$ 
where $\Ker(\Cal C'_1/W)\defeq \Ker (\pi_1(\Cal C'_1)^{\ab}\to \pi_1(W)^{\ab})$, hence an exact sequence 
$$\Ker(\Cal C'_1/W)/l\to \pi_1(\Cal C'_1)^{\ab}/l\to \pi_1(W)^{\ab}/l\to 0.$$ 
By [Katz-Lang], Theorem 1 (together with the fact that $l\in\Primes^{\dag}$), 
$\Ker(\Cal C'_1/W)/l$ is finite, and $\pi_1(W)^{\ab}/l$ is finite by global class field 
theory (cf. [Katz-Lang], Proof of Theorem 4) (resp. as $\pi_1(W)\simeq\widehat{\Bbb Z}$) 
when $p=0$ (resp. $p>0$). Thus, 
$\pi_1(\Cal C'_1)^{\ab}/l$ is finite, hence so is 
$\pi_1(\Cal C')^{\ab}/l$ ($\twoheadleftarrow \pi_1(\Cal C'_1)^{\ab}/l$). 
This finishes the proof of the claim.

By this claim and Hilbert's irreducibility theorem (cf. [Serre2], 10.6), there exists $c\in C^{\cl}$,  
such that the composite map $D_c \hookrightarrow \pi_1(C) \twoheadrightarrow \Pi_{\Cal C}^{A,\Sigma}$, where $D_c\subset \pi_1(C)$
is a decomposition group at $c$ (thus, $D_c$ is only defined up to conjugation and $D_c\isom G_{k(c)}$), is surjective. 
Hence, the natural map 
$H^1(\Pi_{\Cal C}^{A,S},T^{\Sigma}A) \isom 
H^1(\pi_1(C),T^{\Sigma}A)\to H^1(D_c,T^{\Sigma}A)\  (= 
H^1(G_{k(c)},T^{\Sigma}\Cal A_c))$ 
is injective, as desired. 
\qed\enddemo

\definition{Remark 1.3} A similar statement as in Proposition 1.2 holds 
when $T^{\Sigma}A$ is replaced by any finitely generated torsion-free 
$\ZSigma $-module $M$ on which 
$\pi_1(C)$ acts such that the action of $\pi_1(C)$ factors through 
$\pi_1(C)\twoheadrightarrow \pi_1(\Cal C)$ for some model $\Cal C$ of 
$C$ and that 
the weights of the Galois representation 
associated to $M$ 
are distinct from the weight of the cyclotomic character. 
\enddefinition

More generally, we have the following result which generalises Proposition 1.2 to a wider class of base fields.

\proclaim {Proposition 1.4} Let $\Sigma\subset \Primes^{\dag}$ be a 
{\bf finite} subset. 
Assume that $k$ is {\bf Hilbertian} (cf [Serre2], 9.5). 
Then there exists a {\bf finite} subset $S\defeq S(\Sigma)\subset  
C^{\cl}$ of cardinality $\leq 2$, depending on ($A$ and) $\Sigma$, 
such that the natural restriction map
$$
H^1(\pi_1(C),T^{\Sigma}A)\to \prod_{c\in S} H^1(G_{k(c)},T^{\Sigma}\Cal A_c)$$ 
is {\bf injective}. Moreover, in the case where the $\overline k$-trace 
$
\Tr_{K\overline k/\overline k}(A_{K\overline k})$ 
of $A_{K\overline k}\defeq A\times _KK\overline k$ is trivial, there exists such a set $S$ with 
$\sharp (S)=1$. 
\endproclaim


First, we prove the following. 

\proclaim {Proposition 1.5} 
Let $\Sigma\subset \Primes^{\dag}$ be a 
{\bf finite} subset. 
Assume that $k$ is {\bf Hilbertian}. 
Let $M\subset H^1(\pi_1(C),T^{\Sigma}A)$ be a {\bf finitely generated} 
$\ZSigma $-submodule. Then there exists a closed point
$c\in C^{\cl}$, depending on $\Sigma$ and $M$, such that the natural restriction map 
$
M\to H^1(G_{k(c)},T^{\Sigma}\Cal A_c)$ is {\bf injective}.
\endproclaim

\demo{Proof}
Let $L$ be the kernel of the natural map $\pi_1(C)\to\Aut(T^{\Sigma}A)=\prod_{l\in\Sigma}\Aut(T_lA)$. By definition, 
the action of $\pi_1(C)$ on $T^{\Sigma}A$ factors through $\pi_1(C)\twoheadrightarrow \pi_1(C)/L$, hence 
we have the following exact sequence: 
$$
\matrix 0\to H^1(\pi_1(C)/L,T^{\Sigma}A)\to H^1(\pi_1(C),T^{\Sigma}A)\to & H^0(\pi_1(C)/L,H^1(L,T^{\Sigma}A)) \\
&\Vert \\
&\Hom_{\pi_1(C)}(L,T^{\Sigma}A) \\
&\cap \\
&\Hom(L,T^{\Sigma}A).
\endmatrix
$$
Restricting to $M$, we get a homomorphism $M\to \Hom(L,T^{\Sigma}A)$ or, equivalently, a homomorphism 
$L\to\Hom(M,T^{\Sigma}A)$. Let $L_M$ be the kernel of $L\to\Hom(M,T^{\Sigma}A)$. 

\bigskip\noindent
{\bf Claim:} (i) $L_M\subset\pi_1(C)$ is a closed normal subgroup. 

\noindent
(ii) The Frattini subgroup of $\pi_1(C)/L_M$ is open. 

\noindent
(iii) $M\subset H^1(\pi_1(C),T^{\Sigma}A)$ is contained in the image of the 
inflation map 
$$H^1(\pi_1(C)/L_M,T^{\Sigma}A)\hookrightarrow H^1(\pi_1(C),T^{\Sigma}A).$$

Indeed, as $\Hom(L,T^{\Sigma}A)$ is the set of {\it continuous} 
homomorphisms from $L$ to $T^{\Sigma}A$, the subgroup $L_M$ is closed. As the image of 
$M\to \Hom(L,T^{\Sigma}A)$ is contained in $\Hom_{\pi_1(C)}(L,T^{\Sigma}A)$, the subgroup $L_M$ is normal (not only in $L$ but also) 
in $\pi_1(C)$. Thus, (i) follows. Next, we have the following exact sequence of profinite groups: 
$$1\to L/L_M\to \pi_1(C)/L_M\to \pi_1(C)/L\to 1.$$
As 
$$\pi_1(C)/L\hookrightarrow \Aut(T^{\Sigma}A)=\prod_{l\in\Sigma}\Aut(T_lA)\simeq \prod_{l\in\Sigma}\GL_{2d}(\Bbb Z_l)$$
(where $d\defeq\dim(A)$) and 
$$L/L_M\hookrightarrow 
\Hom(M,T^{\Sigma}A)=\prod_{l\in\Sigma}\Hom(M^l,T_lA)\simeq\prod_{l\in\Sigma}\Bbb Z_l^{2dr_l}$$ 
(where 
$M^l\defeq M\otimes_{\widehat{\Bbb Z}^{\Sigma}}\Bbb Z_l$ and 
$r_l\defeq \dim_{\Bbb Q_l}(M^l\otimes_{\Bbb Z_l}\Bbb Q_l)$), (ii) follows from [Serre2], 10.6, Proposition. 
Finally, 
we have the inflation-restriction exact sequence
$$0\to H^1(\pi_1(C)/L_M,T^{\Sigma}A) \to H^1(\pi_1(C),T^{\Sigma}A)\to H^1(L_M,T^{\Sigma}A)$$
arising from the exact sequence $1\to L_M\to \pi_1(C)\to \pi_1(C)/L_M\to 1$. 
By the very definition of $L_M$, the image of $M\subset H^1(\pi_1(C),T^{\Sigma}A)$ 
in 
$H^1(L_M,T^{\Sigma}A)$ 
is trivial. Thus, the assertion of (iii) follows. 
This finishes the proof of the claim. 

As in the proof of Proposition 1.2, 
(ii) of the above claim and the Hilbertian property of $k$ imply (cf. [Serre2], 10.6) that there exists $c\in C^{\cl}$,  
such that the composite map 
$
D_c\ (\isom G_{k(c)})
\hookrightarrow \pi_1(C) \twoheadrightarrow \pi_1(C)/L_M$ is surjective. 
Hence, the composite map 
$H^1(\pi_1(C)/L_M, T^{\Sigma}A)\hookrightarrow H^1(\pi_1(C),T^{\Sigma}A)\to H^1(D_c,T^{\Sigma}A)\  
(= 
H^1(G_{k(c)},T^{\Sigma}\Cal A_c))$ 
is injective. This, together with (iii) of the above claim, finishes the proof of Proposition 1.5. 
\qed\enddemo

\demo{Proof of Proposition 1.4} 
We have the inflation-restriction exact sequence
$$0\to H^1(G_k,T) \to H^1(\pi_1(C),T^{\Sigma}A)\to H^1(\pi_1(C_{\overline k}),T^{\Sigma}A)$$
arising from the exact sequence $1\to \pi_1(C_{\overline k})\to \pi_1(C)\to G_k\to 1$,  
where $T\defeq (T^{\Sigma}A)^{\pi_1(C_{\overline k})}$. 

First, in the special case that $\Tr_{K\overline k/\overline k}(A_{K\overline k})=0$, 
$A(K\overline k)$ is finitely generated by [Lang-N\'eron], 
hence, in particular, $A(K\overline k)^{\tor}$ is finite and 
$T=(T^{\Sigma}A)^{\pi_1(C_{\overline k})}=T^{\Sigma}(A(K\overline k))=0$. 
Thus, the restriction map 
$H^1(\pi_1(C),T^{\Sigma}A)\to H^1(\pi_1(C_{\overline k}),T^{\Sigma}A)$ is 
injective. 
The $\ZSigma $-module 
$H^1(\pi_1(C_{\overline k}),T^{\Sigma}A)
\simeq 
H^1(\Pi_{C_{\overline k}}^{A,\Sigma}, T^{\Sigma}A)$ 
(cf. Lemma 1.1 (iii)) is finitely generated by Lemma 1.6 (ii) below, since $\Pi_{C_{\overline k}}^{A, \Sigma}$ 
is finitely generated as a profinite group (cf. [Grothendieck], Expos\'e XIII, Corollaire 2.12. 
Note that $\Sigma\subset\Primes^{\dag}$). 
As $\ZSigma$ is noetherian, $H^1(\pi_1(C), T^{\Sigma}A)\hookrightarrow 
H^1(\pi_1(C_{\overline k}),T^{\Sigma}A)$ is also finitely generated. 
Thus, the assertion follows from Proposition 1.5 in this case. 

\proclaim{Lemma 1.6} 
{\rm (i)} Let $\Delta$ be a finitely generated group. Let $M$ be a finitely generated 
$\Bbb Z$-module on which $\Delta$ acts. Then $H^1(\Delta, M)$ 
is a finitely generated $\Bbb Z$-module. 

If, moreover, either $\Delta$ or $M$ is finite, then $H^1(\Delta, M)$ is finite. 

\noindent
{\rm (ii)} Let $\Sigma\subset\Primes$ be any subset. 
Let $\Delta$ be a finitely generated profinite group. Let $M$ be a finitely generated 
$\ZSigma $-module on which $\Delta$ acts continuously. Then 
$$H^1(\Delta, M)
\defeq\varprojlim_{\text{$N$:$\Sigma$-integer$>0$}}
H^1(\Delta, M/N)$$ 
is a finitely generated $\ZSigma $-module. 

If, moreover, either $\Delta$ or $M$ is finite, then $H^1(\Delta, M)$ is finite. 
\endproclaim

\demo{Proof}
(i) Take a surjection $(\Bbb Z)^{\oplus s}\twoheadrightarrow M$ of 
$\Bbb Z$-modules and a surjection 
$F_r\twoheadrightarrow \Delta$ of groups, where $F_r=\langle x_1,\dots,x_r\rangle$ 
is a free group of finite rank $r$. Then we claim that 
$H^1(\Delta, M)$ is generated by (at most) $rs$ elements as a 
$\Bbb Z$-module. 
Indeed, 
the inflation map 
$H^1(\Delta, M)\to H^1(F_r,M)$ is injective. By considering the standard resolution of 
the trivial $F_r$-module, we obtain $H^1(F_r,M)=M^{\oplus r}/(x_1-1,\dots,x_r-1)M$, which is 
generated by (at most) $rs$ elements as a $\Bbb Z$-module. As $\Bbb Z$ is 
a PID, $H^1(\Delta, M)\hookrightarrow H^1(F_r,M)$ is also generated by (at most) $rs$ elements 
as a $\Bbb Z$-module, as desired. 

The second assertion follows from the first, together with the standard fact that 
$H^1(\Delta,M)$ is killed by $\sharp(\Delta)$ (resp. $\sharp(M)$) when 
$\sharp(\Delta)<\infty$ (resp. $\sharp(M)<\infty$). 

\noindent 
(ii) Take a surjection $(\ZSigma)^{\oplus s}\twoheadrightarrow M$ of 
$\ZSigma$-modules and a surjection 
$\widehat{F}_r\twoheadrightarrow \Delta$ of profinite groups, where $\widehat{F}_r=\langle x_1,\dots,x_r\rangle$ 
is a free profinite group of finite rank $r$. Then we claim that 
$H^1(\Delta, M)$ is generated by (at most) $rs$ elements as a 
$\ZSigma$-module. Indeed, write $M=\prod_{l\in\Sigma}M_l$ for the canonical 
decomposition corresponding to the decomposition $\ZSigma=\prod_{l\in\Sigma}\Bbb Z_l$. 
Then we have $H^1(\Delta, M)=\prod_{l\in \Sigma}H^1(\Delta, M_l)$, hence 
it suffices to prove that 
$H^1(\Delta, M_l)$ is generated by (at most) $rs$ elements as a 
$\Bbb Z_l$-module.
Now, the inflation map 
$H^1(\Delta, M_l)\to H^1(\widehat{F}_r,M_l)$ is injective. By considering the standard resolution of 
the trivial $\widehat{F}_r$-module, we obtain $H^1(\widehat{F}_r,M_l)=M_l^{\oplus r}/(x_1-1,\dots,x_r-1)M_l$, which is 
generated by (at most) $rs$ elements as a $\Bbb Z_l$-module. As $\Bbb Z_l$ is 
a PID, $H^1(\Delta, M_l)\hookrightarrow H^1(\widehat{F}_r,M_l)$ is also generated by (at most) $rs$ elements 
as a $\Bbb Z_l$-module, as desired. 

The second assertion follows from the first, together with the standard fact that 
$H^1(\Delta,M)$ is killed by $\sharp(\Delta)$ (resp. $\sharp(M)$) when 
$\sharp(\Delta)<\infty$ (resp. $\sharp(M)<\infty$). 
\qed
\enddemo


We shall return to the proof of Proposition 1.4. 
In general, fix any $c_0\in C^{\cl}$. By Proposition 1.5, it suffices to prove that 
the kernel $N$ of 
$
H^1(\pi_1(C),T^{\Sigma}A) \to H^1(G_{k(c_0)},T^{\Sigma}\Cal A_{c_0})$ is finitely generated 
as a $\ZSigma$-module. 
As we have already seen, the $\ZSigma$-module 
$H^1(\pi_1(C_{\overline k}),T^{\Sigma}A)$ is finitely generated, hence, as 
$\ZSigma$ is noetherian, the image of $N$ in 
$H^1(\pi_1(C_{\overline k}),T^{\Sigma}A)$ is also finitely generated. Thus, 
it suffices to prove that the intersection of $N$ and (the image of) $H^1(G_k,T)$  
in 
$H^1(\pi_1(C),T^{\Sigma}A)$ is finitely generated. Here, we have 
$$N\cap H^1(G_k,T) = \Ker(H^1(G_k,T)\to H^1(G_{k(c_0)}, T^{\Sigma}\Cal A_{c_0})).$$ 
Since the map $H^1(G_k,T)\to H^1(G_{k(c_0)}, T^{\Sigma}\Cal A_{c_0})$ 
factors as 
$$H^1(G_k,T)\to H^1(G_{k(c_0)},T)\to H^1(G_{k(c_0)}, T^{\Sigma}\Cal A_{c_0}),$$
it suffices to prove that both 
$N_1\defeq \Ker(H^1(G_k,T)\to H^1(G_{k(c_0)},T))$ 
and 
$N_2\defeq \Ker(H^1(G_{k(c_0)},T)\to H^1(G_{k(c_0)}, T^{\Sigma}\Cal A_{c_0}))$ 
are finitely generated (as $\ZSigma$ is noetherian). 

To prove that $N_1$ is finitely generated, let $k_1/k$ be the normal closure of the 
finite extension $k(c_0)/k$. Then 
$$N_1\subset \Ker(H^1(G_k,T)\to H^1(G_{k_1},T))
\simeq
H^1(\Aut(k_1/k), T^{G_{k_1}}).$$ 
Thus, $N_1$ is finitely generated (in fact, finite) by Lemma 1.6 (ii). 

To prove that $N_2$ is finitely generated, consider the long exact sequence associated 
to the exact sequence $0\to T\to T^{\Sigma}\Cal A_{c_0} \to (T^{\Sigma}\Cal A_{c_0})/T 
\to 0$ of (continuous) $G_{k(c_0)}$-modules. Then we see that there exists a natural surjection 
$$((T^{\Sigma}\Cal A_{c_0})/T )^{G_{k(c_0)}}\twoheadrightarrow N_2,$$ 
{}from which $N_2$ is finitely generated. 
This finishes the proof of Proposition 1.4. 
\qed
\enddemo

\definition {Question 1.7} 
Do the assertions of Proposition 1.2 and Proposition 1.4 also hold 
when $\sharp(\Sigma)=\infty$ (especially, when $\Sigma=\Primes^{\dag}$)? 
\enddefinition

For the cohomology with torsion coefficients, 
one has the following.

%

\proclaim {Proposition 1.8} 
Assume that $k$ is {\bf Hilbertian}. 
Then, for each $\Primes ^{\dag}$-integer $N>0$, 
there exists 
$S\defeq S(N)\subset C^{\cl}$ of cardinality $\leq 2$,
depending on ($A$ and) $N$, such that the natural restriction map 
$$
H^1(\pi_1(C),A[N])\to \prod_{c\in S} H^1(G_{k(c)},\Cal A_c[N])$$ 
is {\bf injective}. 
\endproclaim

First, we prove the following. 

\proclaim {Proposition 1.9} 
Assume that $k$ is {\bf Hilbertian} and let $N$ 
be a $\Primes ^{\dag}$-integer $>0$. 
Let $M\subset H^1(\pi_1(C),A[N])$ be a {\bf finite} 
$\Bbb Z/N$-submodule. Then there exists a closed point
$c\in C^{\cl}$, depending on $N$ and $M$, such that the natural restriction map 
$
M\to H^1(G_{k(c)},\Cal A_c[N])$ is {\bf injective} 
and that the natural restriction map 
$H^0(\pi_1(C),A[N])\to H^0(G_{k(c)},\Cal A_c[N])$
is an isomorphism. 
\endproclaim

\demo{Proof}
The proof is similar to (and even simpler than) that of Proposition 1.5. 
As $A[N]$ is a finite discrete $\pi_1(C)$-module and 
$M\subset H^1(\pi_1(C),A[N])$ is finite, 
it follows from the definition of the profinite group cohomology that 
there exists an open normal subgroup $L_{0,M}\subset \pi_1(C)$ 
such that $L_{0,M}$ acts trivially on $A[N]$ and that 
$M\subset H^1(\pi_1(C),A[N])$ is contained in the image of the 
inflation map 
$$H^1(\pi_1(C)/L_{0,M},A[N])\hookrightarrow H^1(\pi_1(C),A[N]).$$
Now, the Hilbertian property of $k$ implies that there exists $c\in C^{\cl}$,  
such that the composite map $D_c(\isom G_{k(c)})  \hookrightarrow \pi_1(C) \twoheadrightarrow \pi_1(C)/L_{0,M}$ is surjective. 
Hence, the composite map 
$H^1(\pi_1(C)/L_{0,M}, A[N])\hookrightarrow H^1(\pi_1(C),A[N])\to H^1(D_c,A[N])$ 
($= 
H^1(G_{k(c)},\Cal A_c[N])$)  
(resp. $H^0(\pi_1(C)/L_{0,M}, A[N])\isom H^0(\pi_1(C),A[N])\to H^0(D_c,A[N])$ 
($= 
H^0(G_{k(c)},\Cal A_c[N])$)) 
is injective (resp. an isomorphism). 
This finishes the proof of Proposition 1.9. 
\qed\enddemo

\demo{Proof of Proposition 1.8}
The proof is similar to (and even simpler than) that of Proposition 1.4. 
We have the inflation-restriction exact sequence
$$0\to H^1(G_k,T_0) \to H^1(\pi_1(C),A[N])\to H^1(\pi_1(C_{\overline k}),A[N])$$
arising from the exact sequence $1\to \pi_1(C_{\overline k})\to \pi_1(C)\to G_k\to 1$,  
where $T_0\defeq (A[N])^{\pi_1(C_{\overline k})}$. 
Fix any $c_0\in C^{\cl}$. By Proposition 1.9, it suffices to prove that 
the kernel $N_0$ of 
$
H^1(\pi_1(C),A[N]) \to H^1(G_{k(c_0)},\Cal A_{c_0}[N])$ is finite.  
By Lemma 1.6 (ii), together with Lemma 1.1 (iv), 
$H^1(\pi_1(C_{\overline k}),A[N])$ is finite, hence 
the image of $N_0$ in 
$H^1(\pi_1(C_{\overline k}),A[N])$ is also finite. Thus, 
it suffices to prove that the intersection of $N_0$ and 
(the image of) $H^1(G_k,T_0)$ 
in 
$H^1(\pi_1(C),A[N])$ is finite. Here, we have 
$$N_0\cap H^1(G_k,T_0) = \Ker(H^1(G_k,T_0)\to H^1(G_{k(c_0)}, \Cal A_{c_0}[N])).$$ 
Since the map $H^1(G_k,T_0)\to H^1(G_{k(c_0)}, \Cal A_{c_0}[N])$ 
factors as 
$$H^1(G_k,T_0)\to H^1(G_{k(c_0)},T_0)\to 
H^1(G_{k(c_0)}, \Cal A_{c_0}[N]),$$ 
it suffices to prove that both 
$N_{0,1}\defeq \Ker(H^1(G_k,T_0)\to H^1(G_{k(c_0)},T_0))$ 
and 
$N_{0,2}\defeq \Ker(H^1(G_{k(c_0)},T_0)\to H^1(G_{k(c_0)}, \Cal A_{c_0}[N]))$ 
are finite. 

To prove that $N_{0,1}$ is finite, let $k_1/k$ be the normal closure of the 
finite extension $k(c_0)/k$. Then 
$$N_{0,1}\subset \Ker(H^1(G_k,T_0)\to H^1(G_{k_1},T_0))
\simeq 
H^1(\Aut(k_1/k), T_0^{G_{k_1}}).$$ 
Thus, $N_{0,1}$ is finite by Lemma 1.6 (ii). 

To prove that $N_{0,2}$ is finite, consider the long exact sequence associated 
to the exact sequence $0\to T_0\to \Cal A_{c_0}[N] \to \Cal A_{c_0}[N]/T_0 
\to 0$ of 
(discrete) 
$G_{k(c_0)}$-modules. Then we see that there exists a natural surjection 
$$(\Cal A_{c_0}[N]/T_0 )^{G_{k(c_0)}}\twoheadrightarrow N_{0,2},$$ 
{}from which $N_{0,2}$ is finite. 
\qed\enddemo

The following form of N\'eron's specialisation theorem can be obtained 
as an application of the above injectivity results. 

\proclaim {Proposition 1.10} 
Let $N$ be a $\Primes ^{\dag}$-integer $>0$. 

\noindent
{\rm (i)} Assume that $k$ is {\bf Hilbertian}. 
Then there exists 
$S\defeq S(N)\subset C^{\cl}$ of cardinality $\leq 2$,
depending on ($A$ and) $N$, such that the natural specialisation map 
$A(K)/N\to \prod_{c\in S} \Cal A_c(k(c))/N$ 
is {\bf injective}. 

\noindent
{\rm (ii)} Assume that $k$ is {\bf finitely generated} over the prime field 
and infinite. Then there exists a closed point
$c=c(N)\in C^{\cl}$, depending on ($A$ and) $N$, 
such that 
the natural specialisation map $A(K)/N\to\Cal A_c(k(c))/N$ is {\bf injective}, 
that 
the natural specialisation map $A(K)[N]\to\Cal A_c(k(c))[N]$ is an {\bf isomorphism},
and that 
the natural specialisation map $A(K)\to\Cal A_c(k(c))$ is {\bf injective} 
and its cokernel admits no nontrivial $N$-torsion. 
\endproclaim

\demo{Proof} 
For each $c\in C^{\cl}$, we have a natural commutative diagram

$$
\CD
0 @>>> \Cal A(C)/N @>>> H^1(\pi_1(C),A[N]) \\
@. @VVV @VVV \\
0 @>>> \Cal A_c(k(c))/N @>>>  H^1(G_{k(c)},\Cal A_c[N])  
\endCD
$$

\noindent
where the horizontal sequences arise from Kummer exact sequences over (the \'etale site of) $C$ 
and over $k(c)$, and 
the vertical maps are natural specialisation/restriction maps. Note that $\Cal A(C)\isom A(K)$. 
Thus, (i) follows directly from Proposition 1.8. 

Next, if the characteristic $p$ of $k$ is $>0$ 
(resp. $0$), we may replace $C$ by a nonempty open subset 
on which the $p$-rank of fibres of $\Cal A\to C$ 
is constant (resp. by $C$ itself). 
Then, for any $c\in C^{\cl}$, the specialisation map $A(K)^{\tor}\to\Cal A_c(k(c))^{\tor}$ 
is injective. 

Now, assume that $k$ is finitely generated over the prime field, 
then $K$ is also finitely generated over the prime field, hence $A(K)$ is a finitely generated abelian group 
(cf. [Lang-N\'eron]). 
Applying Proposition 1.9 to the finite $\Bbb Z/N$-submodule 
$M=A(K)/N\mosi \Cal A(C)/N \hookrightarrow H^1(G_K,A[N])$, 
we see that there exists $c\in C^{\cl}$ such that 
$A(K)/N\to\Cal A_c(k(c))/N$ is injective and that 
$A(K)[N]\to\Cal A_c(k(c))[N]$ is an isomorphism. 
It follows from these facts, together with [Serre2], 11.1, Criterion, that 
$A(K)\to \Cal A_c(k(c))$ is injective. (Strictly speaking, this argument is applicable 
only when $N\geq 2$. However, the assertion of (ii) for $N=1$ is included in 
that for (any) $N>1$.) Further, as 
$A(K)/N\to\Cal A_c(k(c))/N$ is injective 
and 
$A(K)[N]\to\Cal A_c(k(c))[N]$ is an isomorphism, 
the cokernel of $A(K)\hookrightarrow \Cal A(k(c))$ admits no nontrivial $N$-torsion. 
(Use the Snake Lemma.) 
This finishes the proof of Proposition 1.10. 
\qed\enddemo
        
\subhead
\S2. Selmer Groups 
\endsubhead
We follow the notations in $\S1$. Moreover, 
let $\Sigma\subset \Primes^{\dag}$ be any {\bf nonempty} subset. 

We have a natural commutative diagram

$$
\CD
0 @>>> A(K)^{\wedge,\Sigma} @>>> H^1(G_K,T^{\Sigma}A) @>>> T^{\Sigma}H^1(G_K,A) @>>> 0  \\
@. @VVV      @V
VV    @V
VV \\
0 @>>> \prod_cA_c(K_c)^{\wedge,\Sigma} @>>>  \prod_c H^1(G_{K_c},T^{\Sigma} A_c) @>>> 
\prod_c T^{\Sigma}H^1(G_{K_c},A_c) @>>> 0 \\
\endCD
\tag {2.1}
$$

\noindent
where the horizontal sequences are 
Kummer exact sequences over $K$ and $K_c$, the vertical maps are natural 
restriction maps, and the product is taken over all closed points $c\in C^{\cl}$. 
We have another natural commutative diagram

$$
\CD
0 @>>> \Cal A(C)^{\wedge,\Sigma} @>>> H^1(\pi_1(C),T^{\Sigma}A) @>>> T^{\Sigma}H^1_{\et}(C,\Cal A) @>>> 0  \\
@. @VVV      @V
VV    @V
VV \\
0 @>>> \prod_c \Cal A_c(k(c))^{\wedge,\Sigma} @>>>  \prod_c H^1(G_{k(c)},T^{\Sigma} \Cal A_c) @>>> 
\prod_c T^{\Sigma}H^1(G_{k(c)},\Cal A_c) @>>> 0. \\
\endCD
\tag {2.2}
$$

\noindent
where the horizontal sequences are Kummer exact sequences over (the \'etale site of) $C$ 
and over $k(c)$, 
the vertical maps are natural restriction maps, 
and the product is taken over all closed points $c\in C^{\cl}$. 

\proclaim{Proposition 2.1} 
{\rm (i)} There exists a natural injective map from 
diagram (2.2) to diagram (2.1). 
Further, the maps on the upper left and lower left terms are isomorphisms. 

\noindent
{\rm (ii)} Assume that $k$ is {\bf finitely generated} 
over the prime field. 
Then the map from diagram (2.2) to diagram (2.1) in (i) is an isomorphism. 
\endproclaim

\demo{Proof} 
(i) For the upper rows of diagrams (2.1) and (2.2), we have a natural commutative diagram 
$$
\matrix
0 &\to& \Cal A(C)^{\wedge,\Sigma} &\to& H^1(\pi_1(C),T^{\Sigma}A) &\to& T^{\Sigma}H^1_{\et}(C,\Cal A) 
&\to& 0  \\
&&\downarrow && \downarrow && \downarrow
&& \\
0 &\to& A(K)^{\wedge,\Sigma} &\to& H^1(G_K,T^{\Sigma}A) &\to& T^{\Sigma}H^1(G_K,A) &\to& 0  
\endmatrix
$$
obtained by taking the pullback of \'etale cohomology 
groups via the natural morphism $\Spec(K)\to C$. Here 
the left vertical map coincides with the map induced by 
the natural map $\Cal A(C)\to A(K)$, which is an isomorphism by 
(a variant of) the valuative criterion for properness. The middle 
vertical map is injective by Lemma 1.1 (iii). 
Thus, the right vertical map is also injective. 

For the lower rows of diagrams (2.1) and (2.2), we have a natural commutative diagram 
for each $c\in C^{\cl}$ 
$$
\matrix
0 &\to& \Cal A_c(k(c))^{\wedge,\Sigma} &\to&  H^1(G_{k(c)},T^{\Sigma} \Cal A_c) &\to& 
T^{\Sigma}H^1(G_{k(c)},\Cal A_c) &\to& 0 \\
&&\uparrow && \uparrow && \uparrow && \\
0 &\to& \Cal A_{R_c}(R_c)^{\wedge,\Sigma} &\to&  H^1(\pi_1(R_c),T^{\Sigma} 
A_c
) &\to& 
T^{\Sigma}H^1_{\et}(R_c,\Cal A_{R_c}) &\to& 0 \\
&&\downarrow && \downarrow && \downarrow && \\
0 &\to& A_c(K_c)^{\wedge,\Sigma} &\to& H^1(G_{K_c},T^{\Sigma} A_c) &\to&
T^{\Sigma}H^1(G_{K_c},A_c) &\to& 0 
\endmatrix
\tag{2.3}$$
where $R_c\defeq\widehat{\Cal O}_{C,c}$ (thus, $K_c$ is the field of fractions of $R_c$) 
and 
$\Cal A_{R_c}\defeq\Cal A\times_C R_c$. 
Here, 
the upper vertical maps are 
obtained by taking the pullback of \'etale cohomology 
groups via the natural morphisms $\Spec(k(c))\to\Spec(R_c)$ 
and the modulo $c$ reduction on $\Cal A_{R_c}$, 
and 
the lower vertical maps are 
obtained by taking the pullback of \'etale cohomology 
groups via the natural morphism $\Spec(K_c)\to\Spec(R_c)$. 
The upper left vertical map is an isomorphism, since 
the reduction map $\Cal A_{R_c}(R_c)\to \Cal A_c(k(c))$ is 
surjective with kernel being $N$-divisible for every $\Primes^{\dag}$-integer $N>0$. 
The upper middle vertical map is an isomorphism, as 
$G_{k(c)}\isom \pi_1(R_c)$ and $T^{\Sigma}A_c\isom T^{\Sigma}\Cal A_c$. 
Thus, the upper right vertical map is also an isomorphism. The 
lower left vertical map coincides with the map induced by 
the natural map $\Cal A_{R_c}(R_c)\to A_c(K_c)$, which is an isomorphism by 
the valuative criterion for properness. The middle 
vertical map is injective by Lemma 1.1 (iii). 
Thus, the lower right vertical map is also injective. 

Finally, it follows from various functoriality properties that the above maps form 
a map from diagram (2.2) to diagram (2.1) with the desired properties. 

\noindent
(ii) By Lemma 1.1 (iii), the maps  
$H^1(\pi_1(C),T^{\Sigma}A) \to H^1(G_K,T^{\Sigma}A)$ 
and 
$$H^1(\pi_1(R_c), T^{\Sigma}A_c)(\isom H^1(G_{k(c)},T^{\Sigma} \Cal A_c)) \to 
H^1(G_{K_c},T^{\Sigma} A_c)$$ 
for $c\in C^{\cl}$ are isomorphisms. The assertion follows from this, together with (i). 
\qed\enddemo

\proclaim{Proposition 2.2} Assume that $k$ is {\bf Hilbertian}. 
Then the middle and left vertical maps in diagrams (2.1) and (2.2) are {\bf injective}. 
(
For the kernels of the right vertical maps, see \S 3.) 
\endproclaim

\demo{Proof} 
The middle vertical map in diagram (2.2), say $r_{\Sigma}$, 
is identified with the product of $r_{\{l\}}$ ($l\in\Sigma$). Thus, 
the injectivity of $r_{\Sigma}$ 
follows immediately from Proposition 1.4. 

Next, as in the proof of Lemma 1.1, we have a natural commutative diagram
$$
\matrix
0&\to& H^1(\pi_1(C),T^{\Sigma}A) &\to& H^1(G_K,T^{\Sigma}A) &\to& H^1(I_C,T^{\Sigma}A)  \\
&&\downarrow&&\downarrow&&\downarrow \\
0&\to& \prod_c H^1(G_{k(c)},T^{\Sigma}\Cal A_c )&\to& \prod_c H^1(G_{K_c},T^{\Sigma} A_c) &\to&
\prod_c H^1(I_c,T^{\Sigma} A_c) 
\endmatrix
$$
where the horizontal sequences are inflation-restriction exact sequences, the left and the middle 
vertical maps are the middle vertical maps in diagrams (2.2) and (2.1), respectively, and the 
right vertical map is a natural restriction map. As shown above, the left vertical 
map is injective. The right vertical map is also injective, since 
$H^1(I_C,T^{\Sigma}A)=\Hom(I_C,T^{\Sigma}A)$, 
$H^1(I_c,T^{\Sigma} A_c)=\Hom(I_c,T^{\Sigma} A_c)$, 
$T^{\Sigma} A \isom T^{\Sigma} A_c$, 
and $I_C$ is (topologically) normally generated by $I_c$ ($c\in C^{\cl}$). Thus, 
the middle vertical map (that is, the middle vertical map in diagram (2.1)) is 
injective. 

Finally, the left vertical maps in diagrams (2.1) and (2.2) are injective, as 
the middle vertical maps therein are injective. 
\qed
\enddemo

For the rest of this paper, we will assume that $k$ is {\bf finitely generated} 
over the prime field and infinite. (Thus, in particular, $k$ is Hilbertian.) 
We will identify $A(K)^{\wedge,\Sigma}$, 
$H^1(\pi_1(C),T^{\Sigma}A)$, and $\prod_c\Cal A_c(k(c))^{\wedge,\Sigma}$ with their images in
$\prod_c H^1(G_{k(c)},T^{\Sigma}\Cal A_c)$. 

\definition {Definition 2.3}
(i) For each $\Primes^{\dag}$-integer $N>0$, we define the {\bf $N$-Selmer group} 
$$\Sel_N (A)\defeq \Sel_N (A,C)\defeq \Ker \lgroup H^1(G_K,A[N])\to \prod _c H^1(G_{K_c},A_c)\rgroup.$$ 

\noindent
(ii) We define the {\bf $\Sigma$-adic Selmer group} 
$$
\align
\Sel^{\Sigma} (A)\defeq \Sel^{\Sigma} (A,C)
&\defeq \Ker \lgroup H^1(G_K,T^{\Sigma}A) \to \prod_c T^{\Sigma}H^1(G_{K_c},A_c)\rgroup \\
&=\Ker \lgroup H^1(\pi_1(C),T^{\Sigma}A) \to \prod_c T^{\Sigma}H^1(G_{k(c)},\Cal A_c)\rgroup,
\endalign
$$
so that $\Sel^{\Sigma}(A)=\varprojlim_{\text{$N$:$\Sigma$-integer$>0$}}\Sel_N(A)$. 
\enddefinition

We have natural injective maps 
$A(K)/\{\text{$\Sigma'$-tor}\}\hookrightarrow A(K)^{\wedge,\Sigma}$ 
and 
$\Cal A_c(k(c))/\{\text{$\Sigma'$-tor}\}\hookrightarrow \Cal A_c(k(c))^{\wedge,\Sigma}$ 
($c\in C^{\cl}$), 
as $A(K)$ and $\Cal A_c(k(c))$ are finitely generated $\Bbb Z$-modules (cf. [Lang-N\'eron]). 
We will identify $\prod_c(\Cal A_c(k(c))/\{\text{$\Sigma'$-tor}\})$ 
with its image in $\prod_c \Cal A_c(k(c))^{\wedge,\Sigma}$.

\definition {Definition 2.4} We define the {\bf $\Sigma$-discrete Selmer group}
$$\Se^{\Sigma} (A)\defeq \Se^{\Sigma} (A,C)\defeq 
H^1(\pi_1(C),T^{\Sigma}A)
\bigcap \prod_c \lgroup \Cal A_c(k(c))/\{\text{$\Sigma'$-tor}\}
\rgroup \subset \prod_cH^1(G_{k(c)},T^{\Sigma}A).$$
 \enddefinition

Note that $\Se^{\Sigma}(A)\subset \Sel^{\Sigma} (A)$ by definition.
One of our main results in this section is the following.

\proclaim {Proposition 2.5} 
The $\Sigma$-discrete Selmer group $\Se^{\Sigma}(A)$ is a {\bf finitely generated} $\Bbb Z$-module.
\endproclaim

First, we prove the following.

\proclaim {Lemma 2.6} 
The following holds
$$
\align
&\ H^1(\pi_1(C),T^{\Sigma}A)\bigcap \prod_c(H^1(G_{k(c)},T^{\Sigma}A)^{\tor}) \\
=&\ \Se^{\Sigma}(A)\bigcap \prod_c (\Cal A_c (k(c))^{\tor}/\{\text{$\Sigma'$-$\tor$}\}) \\
=&\ \Se^{\Sigma}(A)^{\tor} 
=A(K)^{\tor}/\{\text{$\Sigma'$-$\tor$}\}.
\endalign$$
\endproclaim

\demo{Proof} 
For each of the three desired equalities, 
the containment relation $\supset$ clearly holds. Thus, it suffices to prove that 
$$H^1(\pi_1(C),T^{\Sigma}A)\bigcap \prod_c(H^1(G_{k(c)},T^{\Sigma}A)^{\tor}) 
\subset
A(K)^{\tor}/\{\text{$\Sigma'$-$\tor$}\}.
$$
So, take any element $\alpha=(\alpha_l)_{l\in \Sigma}$ of 
$H^1(\pi_1(C),T^{\Sigma}A)\bigcap \prod_c(H^1(G_{k(c)},T^{\Sigma}A)^{\tor})$,  
where we write
$\alpha_l\in H^1(\pi_1(C),T_lA)$ for the $l$-component of the cohomology class $\alpha$. 
For each prime $l\in \Sigma$, 
there exists, by Proposition 1.2, 
a 
point $c\in C^{\cl}$ (depending on $l$) 
such that the natural restriction map 
$H^1(\pi_1(C),T_lA)\to H^1(G_{k(c)},T_l\Cal A_c)$ 
is injective. As the injective image of  
$\alpha_l\in H^1(\pi_1(C),T_lA)$ in 
$H^1(G_{k(c)},T_l\Cal A_c)$ lies in 
$H^1(G_{k(c)},T_l\Cal A_c)^{\tor}$, 
we have 
$$\alpha_l\in 
H^1(\pi_1(C),T_lA)^{\tor}=(A(K)^{\wedge,l})^{\tor}
=A(K)^{\tor}/\{\text{$l'$-tor}\}=A(K)^{\tor,l}.$$
Here, the first equality follows from the fact that 
$T_lH^1_{\et}(C,\Cal A)$ is torsion-free, 
the second equality follows from the 
fact that $A(K)$ is a finitely generated $\Bbb Z$-module, 
and the third equality follows as $A(K)^{\tor}$ is a torsion abelian group. 
In particular, $\alpha_l=0$ for all but finitely many $l\in \Sigma$, 
as $A(K)^{\tor}$ is a finite abelian group. 
Now, we conclude that $\alpha = (\alpha_l)_{l\in\Sigma}\in 
(A(K)/\{\text{$\Sigma'$-tor}\})^{\tor}=
A(K)^{\tor}/\{\text{$\Sigma'$-tor}\}$, 
as desired. 
\qed
\enddemo

\proclaim{Proposition 2.7}
Let $\Sigma_1\subset \Sigma_2\subset \Primes^{\dag}$ be nonempty subsets. 
Then there exists a natural exact sequence 
$$0\to A(K)^{\tor,\Sigma_1'}/A(K)^{\tor, \Sigma_2'} \to \Se^{\Sigma_2}(A)\to\Se^{\Sigma_1}(A),$$
where the map $\Se^{\Sigma_2}(A)\to\Se^{\Sigma_1}(A)$ is induced by the projection 
$H^1(\pi_1(C),T^{\Sigma_2}A)\to H^1(\pi_1(C),T^{\Sigma_1}A)$. 
\endproclaim

\demo{Proof} 
We have 
$$
\align
&\ \Ker(\Se^{\Sigma_2}(A)\to\Se^{\Sigma_1}(A)) \\
=&\ \Se^{\Sigma_2}(A)\cap
\Ker\left(
\prod_c \lgroup \Cal A_c(k(c))/\{\text{$\Sigma_2'$-tor}\}\rgroup
\to 
\prod_c \lgroup \Cal A_c(k(c))/\{\text{$\Sigma_1'$-tor}\}\rgroup
\right) \\
=&\ \Se^{\Sigma_2}(A)\cap
\prod_c \lgroup \Cal A_c(k(c))^{\tor, \Sigma_1'}/\Cal A_c(k(c))^{\tor, \Sigma_2'}\rgroup \\
=&\ (A(K)^{\tor}/A(K)^{\tor, \Sigma_2'}) \cap
\prod_c \lgroup \Cal A_c(k(c))^{\tor, \Sigma_1'}/\Cal A_c(k(c))^{\tor, \Sigma_2'}\rgroup \\
=&\ (A(K)^{\tor, \Sigma_1'}/A(K)^{\tor, \Sigma_2'}), 
\endalign
$$
where the third equality follows from Lemma 2.6. 
\qed
\enddemo

The proof of Proposition 2.5 will follow immediately 
from the following, since $\Cal A_c(k(c))$ is a finitely generated $\Bbb Z$-module 
for every $c\in C^{\cl}$ (cf. [Lang-N\'eron]). 

\proclaim {Proposition 2.8} There exists a closed point $ c\in C^{\cl}$ such that the natural restriction map 
$$\Se^{\Sigma}(A)\to \Cal A_{ c}(k( c))/\{\text{$\Sigma'$-$\tor$}\}$$
is {\bf injective}.
\endproclaim

\demo {Proof} 
First, 
let $\Sigma_1\subset \Sigma_2\subset \Primes^{\dag}$ be nonempty subsets, and 
$c\in C^{\cl}$. Then the natural commutative diagram 
$$
\matrix
H^1(\pi_1(C), T^{\Sigma_2}A)&\to& H^1(G_{k(c)}, T^{\Sigma_2}\Cal A_{c}) \\
&&\\
\downarrow && \downarrow \\
&&\\
H^1(\pi_1(C), T^{\Sigma_1}A)&\to& H^1(G_{k(c)}, T^{\Sigma_1}\Cal A_{c}),  
\endmatrix
$$
where the horizontal maps are natural restriction maps and the vertical maps are 
natural projections, 
restricts to a natural commutative diagram
$$
\matrix 
\Se^{\Sigma_2}(A)&\to& \Cal A_{ c}(k( c))/\{\text{$\Sigma_2'$-$\tor$}\}\\
&&\\
\downarrow && \downarrow \\
&&\\
\Se^{\Sigma_1}(A)&\to& \Cal A_{ c}(k( c))/\{\text{$\Sigma_1'$-$\tor$}\}. 
\endmatrix
$$
Now, suppose that the assertion holds for $\Sigma_1$. Then there exists 
a closed point $ c\in C^{\cl}$ such that the lower horizontal map 
of the latter commutative diagram is injective. Thus, 
$$
\align 
&\ \Ker(\Se^{\Sigma_2}(A)\to \Cal A_{ c}(k( c))/\{\text{$\Sigma_2'$-$\tor$}\}) \\
=&\ \Ker(\Se^{\Sigma_2}(A)\to \Cal A_{ c}(k( c))/\{\text{$\Sigma_2'$-$\tor$}\}) 
\cap \Ker(\Se^{\Sigma_2}(A)\to \Se^{\Sigma_1}(A)) \\
=&\ \Ker(\Se^{\Sigma_2}(A)\to \Cal A_{ c}(k( c))/\{\text{$\Sigma_2'$-$\tor$}\}) 
\cap A(K)^{\tor,\Sigma_1'}/A(K)^{\tor, \Sigma_2'} 
= 0, \\
\endalign
$$
where the second equality follows from Proposition 2.7 and the third equality follows from 
the fact that the reduction map $A(K)^{\tor, \Primes^{\dag}}\to\Cal A_{ c}
(k( c))^{\tor, \Primes^{\dag}}$ is injective, as $\Cal A[N]$ is finite \'etale 
over $C$ for any $\Primes^{\dag}$-integer $N>0$. Thus, the assertion also holds for 
$\Sigma_2$. 

Now, to prove the assertion, we may assume that $\Sigma$ is finite, by replacing 
$\Sigma$ with any nonempty finite subset. (We may even assume $\sharp(\Sigma)=1$.) 
Then, by Proposition 1.2, 
there exists a closed point $ c\in C^{\cl}$, 
such that the natural restriction map
$$
H^1(\pi_1(C),T^{\Sigma}A) \to H^1(G_{k( c)},T^{\Sigma}\Cal A_{ c})$$ 
is injective.
In particular, the natural restriction map 
$$\Se^{\Sigma}(A)\to \Cal A_{ c}(k( c))/\{\text{$\Sigma'$-$\tor$}\}$$
is also injective, as desired. \qed
\enddemo


\proclaim {Proposition 2.9} 
{\rm (i)} For each $\Sigma$-integer $N>0$, the natural map 
$\Se^{\Sigma}(A) 
\to H^1(\pi_1(C),A[N])$ 
induces a natural {\bf injective} map 
$$\Se^{\Sigma}(A)/N \hookrightarrow \Sel_N(A)\subset H^1(\pi_1(C),A[N]).$$

\noindent
{\rm (ii)} The natural map $\Se^{\Sigma}(A) 
\to H^1(\pi_1(C),T^{\Sigma}A)$ 
induces a natural {\bf injective} map 
$$\Se^{\Sigma}(A)^{\wedge,\Sigma}\hookrightarrow \Sel^{\Sigma}(A)\subset H^1(\pi_1(C),T^{\Sigma}A).$$ 
\endproclaim

\demo{Proof} 
(i) As $H^1(\pi_1(C),A[N])$ is killed by $N$, 
the natural map $\Se^{\Sigma}(A) 
\to H^1(\pi_1(C),A[N])$ 
induces a natural map $\iota_N: \Se^{\Sigma}(A) /N \to H^1(\pi_1(C),A[N])$. 
As $\Se^{\Sigma}(A)\subset\Sel^{\Sigma}(A)$, the image of $\iota_N$ 
is contained in $\Sel_N(A)$. Thus, it suffices to prove that 
$\iota_N$ is injective. 
Consider the natural exact sequence 
$$0\to \Se^{\Sigma}(A) @>\iota>> H^1(\pi_1(C),T^{\Sigma}A)\to \Coker (\iota) \to 0,$$
where $\iota$ is the natural injection. We claim that $\Coker(\iota)^{\tor,\Sigma}=0$ holds.
Indeed, first we have the equality 
$\Coker (\iota)^{\tor,\Sigma} = ({\Sel^{\Sigma}(A)}/\Se^{\Sigma}(A))^{\tor,\Sigma}$, 
which follows from the fact that $H^1(\pi_1(C),T^{\Sigma}A)/\Sel^{\Sigma} (A)$ injects into 
$\prod_c T^{\Sigma} H^1(G_{k(c)},\Cal A_c)$,  
which is torsion-free. Moreover, 
we have a natural injective map $\Sel^{\Sigma}(A)/\Se^{\Sigma}(A)\to 
\prod_c \Cal A_c(k(c))^{\wedge,\Sigma}/(\Cal A_c (k(c))/\{\text{$\Sigma'$-tor}\})$ 
and the $\Sigma$-torsion of the latter group is trivial, as follows easily from the fact 
that the $\Sigma$-torsion of $\widehat {\Bbb Z}^{\Sigma}/\Bbb Z$ is trivial, and the groups 
$\Cal A_c(k(c))$ are finitely generated. 

Now, for each $\Sigma$-integer $N>0$, we have a commutative diagram of exact sequences
$$
\CD
0 @>>> \Se(A) @>>> H^1(\pi_1(C),T^{\Sigma}A) @>>> \Coker (\iota) @>>> 0\\
@. @VNVV   @VNVV   @VNVV \\
0 @>>> \Se(A) @>>> H^1(\pi_1(C),T^{\Sigma}A) @>>> \Coker (\iota) @>>> 0\\
\endCD
$$
where the vertical maps are the maps of multiplication by $N$. 
Thus, by the Snake Lemma, we have a natural exact sequence
$$0=\Coker(\iota)[N] 
\to \Se^{\Sigma}(A)/N
\to H^1(\pi_1(C),T^{\Sigma}A)/N.$$ 
Now the assertion follows, as $H^1(\pi_1(C),T^{\Sigma}A)/N\hookrightarrow H^1(\pi_1(C),A[N])$. 

\noindent
(ii) By (i), the natural map 
$$
\matrix
\Se^{\Sigma}(A)^{\wedge,\Sigma}&&
\Sel^{\Sigma}(A)
&&H^1(\pi_1(C),T^{\Sigma}A)\\
\Vert&&\Vert&&\Vert\\
\varprojlim_N
\Se^{\Sigma}(A)/N 
&\to&
\varprojlim_N
\Sel_N(A)
&\subset&
\varprojlim_N
H^1(\pi_1(C),A[N]), 
\endmatrix
$$
where $N$ runs over all $\Sigma$-integers $>0$, is injective, as desired. 
\qed
\enddemo

\proclaim {Proposition 2.10} 
For each $\Primes^{\dag}$-integer $N>0$, the $N$-Selmer group
$\Sel_N(A)$ is {\bf finite}. 
\endproclaim

Let $N$ be a $\Primes ^{\dag}$-integer $>0$. 
As $\Cal A[N]$ is finite \'etale over $C$, 
the $G_K$-module $A[N]$ is unramified, i.e.,
the group $I_C=\Ker(G_K\twoheadrightarrow\pi_1(C))$ acts trivially on $A[N]$. Thus, $A[N]$ has a natural structure of $\pi_1(C)$-module,
and we have a natural inflation-restriction exact sequence
$$0\to H^1(\pi_1(C),A[N]) @>\inf>>   H^1(G_K,A[N]) @> \res>> \Hom (I_C,A[N])^{\pi_1(C)}.$$

\proclaim {Lemma 2.11} The following holds: $\Sel_N(A)\subset H^1(\pi_1(C),A[N])$.
\endproclaim

\demo{Proof} 
As in diagram (2.3), we have a natural commutative diagram 
for each $c\in C^{\cl}$ 
$$
\matrix
0 &\to& \Cal A_{R_c}(R_c)/N &\to&  H^1(\pi_1(R_c),\Cal A_{R_c}[N]) &\to& 
H^1_{\et}(R_c,\Cal A_{R_c})[N] &\to& 0 \\
&&\downarrow && \downarrow && \downarrow && \\
0 &\to& A_c(K_c)/N &\to& H^1(G_{K_c}, A_c[N]) &\to&
H^1(G_{K_c},A_c)[N] &\to& 0, 
\endmatrix
$$
where 
the 
vertical maps are 
obtained by taking the pullback of \'etale cohomology 
groups via the natural morphism $\Spec(K_c)\to\Spec(R_c)$. 
The 
left vertical map coincides with the map induced by 
the natural map $\Cal A_{R_c}(R_c)\to A_c(K_c)$, which is an isomorphism by 
the valuative criterion for properness. 
The 
middle 
vertical map is injective by Lemma 1.1 (iv). 
Thus, the 
right vertical map is also injective. 
By definition, the image of 
$\Sel_N(A)$ in $H^1(G_{K_c},A[N])$ is contained in 
(the image of) $A_c(K_c)/N\mosi \Cal A_{R_c}(R_c)/N$, 
hence, in particular, in (the image of) 
$H^1(\pi_1(R_c),\Cal A_{R_c}[N])$. 
It follows from this that the image of $\Sel_N(A)$ in $ H^1(I_c,A[N])=\Hom(I_c,A[N])$ 
is trivial for all $c\in C^{\cl}$.
Thus, the image of $\Sel_N(A)$ in  $\Hom (I_C,A[N])$ is trivial, as desired, 
since $I_C$ is (topologically) normally generated by the various $I_c$ for $c\in C^{\cl}$.
\qed\enddemo

\demo {Proof of Proposition 2.10} By Proposition 1.8, 
there exists 
$S\subset C^{\cl}$ of cardinality $\leq 2$,
depending on ($A$ and) $N$, such that the natural restriction map 
$$
H^1(\pi_1(C),A[N])\to \prod_{c\in S} H^1(G_{k(c)},\Cal A_c[N])$$ 
is injective. 
In particular, $\Sel_N(A)$ injects into 
$\prod_{c\in S} \Cal A_c(k(c))/N$, which is finite as $S$ is finite and 
$\Cal A_c(k(c))$ is finitely generated for each $c$. This finishes the proof of 
Proposition 2.10.
\qed
\enddemo

\subhead
\S3. The Shafarevich-Tate Group
\endsubhead
We use the same notations as in $\S1$ and $\S2$. 
In particular, 
$k$ is a field which is {\bf finitely generated} over the prime field of 
characteristic $p\geq 0$ and infinite; 
$C\to \Spec k$ is a smooth, separated and geometrically connected 
algebraic curve over $k$ with 
function field $K=k(C)$; and 
$\Cal A\to C$ is an abelian scheme over $C$ with 
generic fibre $A= \Cal A_K = \Cal A \times _C \Spec K$. 

\definition {Definition 3.1} We define the {\bf Shafarevich-Tate group}
$$\Sha(A)\defeq \Sha (A,C)\defeq \Ker \lgroup H^1(G_K,A)\to \prod_c H^1(G_{K_c},A_c)\rgroup,$$
where the product is taken over all closed points $c\in C^{\cl}$.
We set $\Sha(A)^{(\dag)}\defeq \Sha(A)$ 
(resp. $\Sha(A)^{(\dag)}\defeq \Sha(A)/\{\text{$p$-tor}\}$), 
when the characteristic $p$ of $k$ is $0$ (resp. $>0$). 
\enddefinition

Note that the abelian group $\Sha(A)$ is a torsion group since the Galois cohomology group 
$H^1(G_K,A)$ is torsion. 
(In particular, $\Sha(A)^{(\dag)}$ is naturally identified with 
$\Sha(A)^{\tor,p'}$ ($\subset\Sha(A)$) when $p>0$.) 
For each $\Primes^{\dag}$-integer $N>0$, we have a natural exact sequence 
$$0 \to A(K)/N \to \Sel_N(A) \to \Sha (A)[N] \to 0,$$
and 
$$0 \to A(K)^{\wedge,\Sigma} \to \Sel^{\Sigma}(A) \to T^{\Sigma}\Sha (A)\to 0.$$

\definition {Definition 3.2} We define the {\bf $\Sigma$-discrete Shafarevich-Tate group}
$$\Sh^{\Sigma}(A)\defeq \Sh^{\Sigma}(A,C)\defeq \Se^{\Sigma}(A)/(A(K)/\{\text{$\Sigma'$-tor}\}).$$
We set $\Sh^{\Sigma}(A)^{(\dag)}\defeq \Sh^{\Sigma}(A)$ 
(resp. $\Sh^{\Sigma}(A)^{(\dag)}\defeq \Sh^{\Sigma}(A)/\{\text{$p$-tor}\}$), 
when the characteristic $p$ of $k$ is $0$ (resp. $>0$). 
\enddefinition

By definition, we have a natural exact sequence 
$$0 \to A(K)/\{\text{$\Sigma'$-tor}\} \to \Se^{\Sigma}(A) \to \Sh^{\Sigma} (A)\to 0.$$

\proclaim{Proposition 3.3} 
The $\Sigma$-discrete Shafarevich-Tate group $\Sh^{\Sigma}(A)$ is a {\bf finitely generated} 
$\Bbb Z$-module. Further, $\Sh^{\Sigma}(A)^{(\dag)}$ is a {\bf finitely generated free}  $\Bbb Z$-module. 
\endproclaim

\demo{Proof}
The first assertion follows immediately from Proposition 2.5. To prove the second assertion, 
it suffices to show that for each prime number $l\in\Primes^{\dag}$, 
$\Sh^{\Sigma}(A)$ admits no nontrivial $l$-torsion. It follows from Lemma 2.6, together with 
the Snake Lemma, that this last condition is equivalent to the injectivity of the natural map 
$(A(K)/\{\text{$\Sigma'$-tor}\})/l \to \Se^{\Sigma}(A)/l$. By definition, we have 
a natural map $\Se^{\Sigma}(A)\to \Cal A_c(k(c))/\{\text{$\Sigma'$-tor}\}$ for each $c\in C^{\cl}$ 
whose composite with the natural map $A(K)/\{\text{$\Sigma'$-tor}\} \to \Se^{\Sigma}(A)$ 
coincides with the specialisation map 
$A(K)/\{\text{$\Sigma'$-tor}\}=\Cal A(C)/\{\text{$\Sigma'$-tor}\} \to 
\Cal A_c(k(c))/\{\text{$\Sigma'$-tor}\}$. Thus, to prove the desired injectivity, it suffices to 
show that for each prime number $l\in\Primes^{\dag}$, there exists $c\in C^{\cl}$ 
(which may depend on $l$), such that 
the specialisation map 
$(A(K)/\{\text{$\Sigma'$-tor}\})/l
\to 
(\Cal A_c(k(c))/{(\text{$\Sigma'$-tor})})/l$ 
is injective. This last assertion 
follows from Proposition 1.10 (ii). Indeed, by Proposition 1.10 (ii), 
there exists a closed point
$c\in C^{\cl}$, 
such 
that 
the natural specialisation map $A(K)/l\to\Cal A_c(k(c))/l$ is injective, 
and 
that 
the natural specialisation map $A(K)\to\Cal A_c(k(c))$ is injective 
and the cokernel 
$\Cal A_c(k(c))/A(K)$ admits no nontrivial $l$-torsion. 
Further, as $\Cal A_c(k(c))^{\tor,\Sigma'}/A(K)^{\tor,\Sigma'} \hookrightarrow \Cal A_c(k(c))/A(K)$, 
the cokernel $\Cal A_c(k(c))^{\tor,\Sigma'}/A(K)^{\tor,\Sigma'}$ also admits no nontrivial $l$-torsion, 
or, equivalently, we have 
$A(K)^{\tor,\Sigma'}[l^{\infty}]\isom \Cal A_c(k(c))^{\tor,\Sigma'}[l^{\infty}]$. 
Consider the following 
commutative diagram of exact sequences: 
$$
\CD
0 @>>> A(K)^{\tor,\Sigma'}/l @>>> A(K)/l @>>> (A(K)/\{\text{$\Sigma'$-tor}\})/l @>>> 0 \\
@. @VVV @VVV @VVV @. \\
0 @>>> \Cal A_c(k(c))^{\tor,\Sigma'}/l @>>> \Cal A_c(k(c))/l @>>> (\Cal A_c(k(c))/\{\text{$\Sigma'$-tor}\})/l @>>> 0 
\endCD
$$
where the vertical maps are natural specialisation maps and the injectivity of 
the horizontal map 
$A(K)^{\tor,\Sigma'}/l \to A(K)/l$ 
(resp. $\Cal A_c(k(c))^{\tor,\Sigma'}/l \to \Cal A_c(k(c))/l$) follows from the 
fact that 
$A(K)^{\tor,\Sigma'}/l=0$ (resp. $\Cal A_c(k(c))^{\tor,\Sigma'}/l=0$) 
if $l\in \Sigma$ and 
$(A(K)/\{\text{$\Sigma'$-tor}\})[l]=0$ (resp. 
$(\Cal A_c(k(c))/\{\text{$\Sigma'$-tor}\})[l]=0$) if $l\in\Sigma'$. 
Thus, by the Snake Lemma, to prove the injectivity of the map 
$(A(K)/\{\text{$\Sigma'$-tor}\})/l
\to 
(\Cal A_c(k(c))/\{\text{$\Sigma'$-tor}\})/l$, 
it suffices to prove that the map 
$A(K)^{\tor,\Sigma'}/l\to \Cal A_c(k(c))^{\tor,\Sigma'}/l$ 
is an isomorphism, which follows from the above-mentioned isomorphism 
$A(K)^{\tor,\Sigma'}[l^{\infty}]\isom \Cal A_c(k(c))^{\tor,\Sigma'}[l^{\infty}]$. 
This finishes the proof of Proposition 3.3. 
\qed
\enddemo

\proclaim{Proposition 3.4} 
Let $\Sigma_1\subset \Sigma_2\subset \Primes^{\dag}$ be nonempty subsets. 
Then the natural projection  $T^{\Sigma_2}\Sha(A) \twoheadrightarrow T^{\Sigma_1}\Sha(A)$ 
induces an {\bf injective} map $\Sh^{\Sigma_2}(A)\hookrightarrow \Sh^{\Sigma_1}(A)$. 
\endproclaim

\demo{Proof} 
Consider the following 
commutative diagram of exact sequences: 
$$
\CD
0 @>>> A(K)/\{\text{$\Sigma_2'$-tor}\} @>>> \Se^{\Sigma_2}(A) @>>> \Sh^{\Sigma_2}(A) @>>> 0 \\
@. @VVV @VVV @VVV @. \\
0 @>>> A(K)/\{\text{$\Sigma_1'$-tor}\} @>>> \Se^{\Sigma_1}(A) @>>> \Sh^{\Sigma_1}(A) @>>> 0 \\
\endCD
$$
where the left vertical map is (resp. the middle and right vertical maps are)  
induced by the identity map $A(K)\to A(K)$ (resp. the natural 
projection 
$H^1(\pi_1(C),T^{\Sigma_2}A)\to H^1(\pi_1(C),T^{\Sigma_1}A)$). 
By the Snake Lemma, the desired injectivity follows from Proposition 2.7, together with the 
surjectivity of the left vertical map 
$A(K)/\{\text{$\Sigma_2'$-tor}\}\to A(K)/\{\text{$\Sigma_1'$-tor}\}$. 
\qed\enddemo

\proclaim{Proposition 3.5} 
{\rm (i)} For each $\Sigma$-integer $N>0$, the natural map 
$\Sh^{\Sigma}(A) 
\to \Sha(A)[N]$ 
induces a natural {\bf injective} map 
$$\Sh^{\Sigma}(A)/N \hookrightarrow \Sha(A)[N]\subset H^1_{\et}(C,\Cal A)[N].$$

\noindent
{\rm (ii)} The natural map $\Sh^{\Sigma}(A) 
\to T^{\Sigma}\Sha(A)$ 
induces a natural {\bf injective} map 
$$\Sh^{\Sigma}(A)^{\wedge,\Sigma}\hookrightarrow T^{\Sigma}\Sha(A)\subset T^{\Sigma}H^1_{\et}(C,\Cal A).$$ 
\endproclaim

\demo{Proof}
(i) Consider the following commutative diagram of exact sequences: 
$$
\CD
 @. A(K)/N @>>> \Se^{\Sigma}(A)/N @>>> \Sh^{\Sigma}(A)/N @>>> 0 \\
@. @VVV @VVV @VVV @. \\
0 @>>> A(K)/N @>>> \Sel_N(A) @>>> \Sha(A)[N] @>>> 0 \\
\endCD
$$
where the left (resp. middle, resp. right) vertical map is the identity map 
(resp. induced by  the natural map 
$\Se^{\Sigma}(A) \to H^1(\pi_1(C),A[N])$ as in Proposition 2.9 (i), 
resp. induced by the natural map $\Sh^{\Sigma}(A) \to \Sha(A)[N]$). It follows 
from this first that
the left upper horizontal map $A(K)/N \to \Se^{\Sigma}(A)/N$ is injective. 
Now, applying the Snake Lemma to the above diagram, 
we see that the desired injectivity follows from 
Proposition 2.9 (i). 

\noindent
(ii) 
By (i), the natural map 
$$
\matrix
\Sh^{\Sigma}(A)^{\wedge,\Sigma}&&
T^{\Sigma}\Sha(A)
&&T^{\Sigma} H^1(C,\Cal A)\\
\Vert&&\Vert&&\Vert\\
\varprojlim_N
\Sh^{\Sigma}(A)/N 
&\to&
\varprojlim_N
\Sha(A)[N]
&\subset&
\varprojlim_N
H^1(C,\Cal A)[N], 
\endmatrix
$$
where $N$ runs over all $\Sigma$-integers $>0$, is injective, as desired. 
\qed
\enddemo

%
%

\proclaim{Proposition 3.6}
The natural map $\Sh^{\Sigma}(A) \to T^{\Sigma}\Sha(A)$ 
induces a natural {\bf injective} map 
$$\Sh^{\Sigma}(A)^{(\dag)}\hookrightarrow T^{\Sigma}\Sha(A)\subset T^{\Sigma}H^1_{\et}(C,\Cal A).$$ 
\endproclaim

\demo{Proof}
Consider the following commutative diagram of exact sequences: 
$$\CD
0 @>>> \Sh^{\Sigma}(A)^{\tor,\dag'} @>>> \Sh^{\Sigma}(A) @>>> \Sh^{\Sigma}(A)^{(\dag)} @>>> 0 \\
@. @VVV @VVV @VVV @. \\
 @. (\Sh^{\Sigma}(A)^{\tor,\dag})^{\wedge,\Sigma} @>>> \Sh^{\Sigma}(A)^{\wedge,\Sigma} @>>> (\Sh^{\Sigma}(A)^{(\dag)})^{\wedge,\Sigma} @>>> 0 
\endCD
$$
where $\Sh^{\Sigma}(A)^{\tor,\dag'}$ stands for $\Sh^{\Sigma}(A)^{\tor,(\Primes^{\dag})'}$ 
(i.e., 
$\Sh^{\Sigma}(A)^{\tor,\dag'}=\Sh^{\Sigma}(A)^{\tor,p}$ for $p>0$ and 
$\Sh^{\Sigma}(A)^{\tor,\dag'}=0$ for $p=0$) and the 
lower horizontal sequence is the $\Sigma$-adic completion of the upper horizontal sequence. 
As $\Sigma\subset\Primes^{\dag}$, we have $(\Sh^{\Sigma}(A)^{\tor,\dag'})^{\wedge,\Sigma}=0$, hence 
$\Sh^{\Sigma}(A)^{\wedge,\Sigma} \isom (\Sh^{\Sigma}(A)^{(\dag)})^{\wedge,\Sigma}$. 
Thus, the desired injectivity is equivalent (cf. Proposition 3.5 (ii)) 
to the injectivity of the right vertical map 
$\Sh^{\Sigma}(A)^{(\dag)}\to (\Sh^{\Sigma}(A)^{(\dag)})^{\wedge,\Sigma}$, which follows from 
the fact (cf. Proposition 3.3) that 
$\Sh^{\Sigma}(A)^{(\dag)}$ is a finitely generated free $\Bbb Z$-module, as $\Sigma\nemp$. 
\qed\enddemo

\proclaim {Proposition 3.7}
Assume that there exists $l\in\Sigma$ such that $T_l\Sha(A)=0$. Then $\Sh^{\Sigma}(A)^{(\dag)}=0$.
\endproclaim

\demo{Proof} 
We have 
$$\Sh^{\Sigma}(A)^{(\dag)}\hookrightarrow \Sh^{\{l\}}(A)^{(\dag)}
\hookrightarrow T_l\Sha(A)=0,$$
by Proposition 3.4 
and Proposition 3.6. Thus, the assertion follows. 
\qed
\enddemo

We conjecture the following.

\definition {Conjecture 3.8} We have 
$\Sh^{\Sigma}(A)=0$ or, equivalently, $\Se^{\Sigma}(A)=A(K)/\{\text{$\Sigma'$-$\tor$}\}$ 
unconditionally.  
\enddefinition


\proclaim {Proposition 3.9} 
Let $N$ be a $\Primes^{\dag}$-integer $>0$. Then: 

\noindent
{\rm (i)} 
$\Sha(A)[N]$ is {\bf finite}.  

\noindent
{\rm (ii)} $\Sha(A)/N$ is {\bf finite}. 
\endproclaim

\demo{Proof}
(i) This follows from Proposition 2.10, 
as we have an exact sequence 
$$0 \to A(K)/N \to \Sel_N(A) \to \Sha (A)[N] \to 0.$$

\noindent
(ii) This follows from (i), together with Lemma 3.10 below. 
\qed\enddemo

\proclaim {Lemma 3.10} 
Let $M$ be a torsion abelian group and $N>0$ an integer. 
Assume that $M[N]$ is finite. Then $M/N$ is finite. 
\endproclaim

\demo{Proof} 
There are several ways of proving this elementary fact. 
For example, consider the following decreasing sequence: 
$$M[N]\supset NM[N^2]\supset \cdots\supset N^{n-1}M[N^n]\supset N^nM[N^{n+1}]\supset\cdots,$$ 
which stabilises as $M[N]$ is finite. Thus, 
there exists an integer $n_0>0$, such that $C\defeq N^{n_0-1}M[N^{n_0}]=N^{n-1}M[N^n]$ $(\subset M[N])$ 
for all $n> n_0$. Now, for each $n> n_0$, consider the following commutative diagram of 
exact sequences: 
$$
\CD
0 @>>> M[N^{n-1}] @>>> M[N^n] @>{N^{n-1}}>> C @>>> 0 \\
@. @V{N}VV @V{N}VV @V{N=0}VV @. \\
0 @>>> M[N^{n-1}] @>>> M[N^n] @>{N^{n-1}}>> C @>>> 0. \\
\endCD
$$
By the Snake Lemma, we have an exact sequence: 
$$0\to C \to M[N^{n-1}]/N\to M[N^n]/N\to C\to 0.$$
Thus, in particular, $\sharp(M[N^n]/N)$ stabilises and is bounded, hence 
$$M/N=M[N^{\infty}]/N=(\varinjlim_{n\geq 0} M[N^n])/N=\varinjlim_{n\geq 0} (M[N^n]/N)$$ 
is finite, as desired. 
\qed\enddemo

\proclaim {Proposition 3.11} 
Let $A\to A'$ be an isogeny of abelian varieties over $K$. Then it induces 
a natural homomorphism $\Sha(A)=\Sha(A,C)\to\Sha(A',C)=\Sha(A')$, hence 
a natural homomorphism $\Sha(A)^{(\dag)}\to\Sha(A')^{(\dag)}$. 
Further, the kernel and the cokernel of the latter homomorphism 
are both finite. 
\endproclaim

\demo{Proof}
By definition, any homomorphism $A\to A'$ over $K$ induces $\Sha(A)\to\Sha(A')$, 
hence $\Sha(A)^{(\dag)}\to\Sha(A')^{(\dag)}$, functorially. 
Now, if $f:A\to A'$ is an isogeny, then there exist an isogeny $g:A'\to A$ and an integer $N>0$, 
such that $g\circ f=N\cdot \id_A$ and $f\circ g=N\cdot \id_{A'}$. 
We define $N^{\dag}$ to be the maximal $\Primes^{\dag}$-integer dividing $N$. 
(Thus, $N/N^{\dag}$ is $1$ (resp. the maximal $p$-power dividing $N$) 
if $p=0$ (resp. $p>0$).) 
By functoriality, these equalities 
imply $\Ker(\Sha(A)\to\Sha(A'))\subset\Sha(A)[N]$, hence 
$\Ker(\Sha(A)^{(\dag)}\to\Sha(A')^{(\dag)})\hookrightarrow \Sha(A)[N^{\dag}]$, and 
$\Sha(A')/N\twoheadrightarrow \Coker(\Sha(A)\to\Sha(A'))$, hence 
$\Sha(A')/N^{\dag}\twoheadrightarrow \Coker(\Sha(A)^{(\dag)}\to\Sha(A')^{(\dag)})$. 
Now, the desired finiteness follows from Proposition 3.9. 
\qed\enddemo

\proclaim {Proposition 3.12} 
Let $k'/k$ be a finite extension of fields (finitely generated over the prime field), 
$C'\to \Spec k'$ a smooth, separated and geometrically connected algebraic curve over $k'$, 
and $C'\to C$ a dominant $k$-morphism. Write $K'=k'(C')$ for the function field of $C'$, 
and let $K'/K$ be the (finite) extension of function fields induced by $C'\to C$. 
Then it induces a natural homomorphism $\Sha(A)=\Sha(A,C)\to\Sha(A_{K'},C')=\Sha(A_{K'})$, 
hence a natural homomorphism $\Sha(A)^{(\dag)}\to\Sha(A_{K'})^{(\dag)}$. Further, the kernel of 
the former (resp. latter) homomorphism is finite, if $K'/K$ is separable (resp. in general). 
\endproclaim

\demo{Proof}
By definition, $C'\to C$ induces 
$\Sha(A)=\Sha(A,C)\to\Sha(A_{K'},C')=\Sha(A_{K'})$, 
hence $\Sha(A)^{(\dag)}\to\Sha(A_{K'})^{(\dag)}$, functorially. 
First, assume that $K'/K$ is separable. Then, replacing 
$K'/K$ by its Galois closure $K''/K$, 
$k'$ by the algebraic closure $k''$ of $k$ in $K''$ and 
$C'$ by the smooth locus over $k''$ of the integral closure of $C'$ in $K''$, 
if necessary, we may reduce 
the finiteness of $\Ker(\Sha(A)\to\Sha(A_{K'}))$ to the case 
where $K'/K$ is Galois. In this case, we have 
$$\Ker(\Sha(A)\to\Sha(A_{K'}))\subset
\Ker(H^1(G_K,A)\to H^1(G_{K'},A_{K'}))=
H^1(\Gal(K'/K), A(K')),$$
which is finite by Lemma 1.6 (i) (together with [Lang-N\'eron]), as desired. 

Thus, to prove the finiteness of 
$\Ker(\Sha(A)^{(\dag)}\to\Sha(A_{K'})^{(\dag)})$ 
in general, we may assume $p>0$. Then, for any finite 
extension $K'/K$, there exist $n\geq 0$ and a finite 
separable extension $K''/K^{\frac{1}{p^n}}$ such that 
$K''\supset K'$. The algebraic closures of $k$ in 
$K'$ and $K^{\frac{1}{p^n}}$ are $k'$ and $k^{\frac{1}{p^{n}}}$, 
respectively, and let $k''$ be the algebraic closure of 
$k$ in $K''$, which is separable over $k^{\frac{1}{p^n}}$. 
Let $\tilde C'$, $C^{\frac{1}{p^n}}$ and 
$\tilde C''$ be the integral closure of $C$ in $K'$, 
$K^{\frac{1}{p^n}}$ and $K''$, respectively. 
Then $\tilde C'$, $C^{\frac{1}{p^n}}$ and $\tilde C''$ 
are regular, separated, geometrically connected curves 
over $k'$, $k^{\frac{1}{p^n}}$ and $k''$, respectively. 
Further, $\tilde C'\supset C'$ is generically smooth over $k'$, 
$(C^{\frac{1}{p^n}}\to\Spec k^{\frac{1}{p^n}})
\simeq (C\to\Spec k)$ is smooth, and 
$\tilde C''$ is generically \'etale over $C^{\frac{1}{p^n}}$, hence 
generically smooth over $k^{\frac{1}{p^n}}$, and over $k''$. 
Let $C''\subset \tilde C''$ be the intersection of the smooth 
locus of $\tilde C''\to \Spec k''$ and the inverse image of 
$C'$ under the (finite) morphism $\tilde C''\to \tilde C'$. 
Then $C''$ is smooth, separated, geometrically connected curve 
over $k''$. 

Now, the natural homomorphism $\Sha(A,C)\to\Sha(A_{K''},C'')$ 
factors in two ways, as 
$\Sha(A,C)\to\Sha(A_{K'},C')\to \Sha(A_{K''},C'')$ 
and as 
$\Sha(A,C)
=\Sha(A_{K^{\frac{1}{p^0}}},C^{\frac{1}{p^0}})
\to \Sha(A_{K^{\frac{1}{p}}},C^{\frac{1}{p}})
\to \Sha(A_{K^{\frac{1}{p^2}}},C^{\frac{1}{p^2}})
\to\cdots
\to \Sha(A_{K^{\frac{1}{p^n}}},C^{\frac{1}{p^n}})
\to \Sha(A_{K''},C'')$. It follows from this that 
the finiteness of 
$\Ker(\Sha(A)^{(\dag)}\to\Sha(A_{K'})^{(\dag)})$ 
is reduced to that of 
$\Ker(\Sha(A_{K^{\frac{1}{p^i}}},C^{\frac{1}{p^{i}}})^{(\dag)}
\to \Sha(A_{K^{\frac{1}{p^{i+1}}}},C^{\frac{1}{p^{i+1}}})^{(\dag)})$ 
($i=0,1,\dots,n-1$) 
and that of 
$\Ker(\Sha(A_{K^{\frac{1}{p^n}}},C^{\frac{1}{p^n}})^{(\dag)}
\to \Sha(A_{K''},C'')^{(\dag)})$. The latter finiteness 
follows from the above argument, 
as $K''/K^{\frac{1}{p^n}}$ is separable. 
For the former finiteness, it suffices to prove it 
for $i=1$. 
In this case, the inclusion $K\subset K^{\frac{1}{p}}$ can be identified 
with the inclusion $\sigma:K\hookrightarrow K$, $x\mapsto x^p$. 
Under this identification, we have $A_{K^{\frac{1}{p}}}=A\times_{K,\sigma}K$ and 
the homomorphism $\Sha(A,C)\to\Sha(A_{K^{\frac{1}{p}}},C^{\frac{1}{p}})$ is identified with the 
homomorphism $\Sha(A,C)\to\Sha(A\times_{K,\sigma}K,C)$ induced by the relative Frobenius $K$-morphism 
$A\to A\times_{K,\sigma}K$. Thus, the desired finiteness of 
$\Ker(\Sha(A)^{(\dag)}\to\Sha(A_{K^{\frac{1}{p}}})^{(\dag)})$ follows from Proposition 3.11. 
\qed\enddemo

\definition{Remark 3.13} 
For simplicity, write $Q$ for the prime field of $k$ and $Z$ for the image of $\Bbb Z$ in $Q$. 
Thus, we have $Q=\Bbb Q$ (resp. $Q=\Bbb F_p$) and $Z=\Bbb Z$ (resp. $Z=\Bbb F_p$) when 
$p=0$ (resp. $p>0$). Then, as $k$ is finitely generated over the perfect field $Q$, the system 
$
C \to \Spec k\to\Spec Q$ admits a smooth model 
$
\Cal C \to V \to U$. More precisely, $U=\Spec Z$; $V$ is an integral scheme which is smooth over $U$ 
and whose function field is isomorphic to (and is identified with) $k$; $\Cal C$ is a smooth scheme over $V$ whose 
generic fibre $\Cal C\times_V k$ is $k$-isomorphic to (and is identified with) $C$. 
Let $\Cal C^1$ denote the set of points of codimension $1$ of $\Cal C$, hence we have $C^{\cl}\subset\Cal C^1$. 
For each $c\in \Cal C^1$, 
let 
$K_c$ be the completion of $K$ at $c$, and $A_c\defeq A\times _KK_c$, just as in the case of $c\in C^{\cl}$. 
We define 
$$\Sel_N (A,\Cal C)\defeq \Ker \lgroup H^1(G_K,A[N])\to \prod _{c\in\Cal C^1} H^1(G_{K_c},A_c)\rgroup$$ 
for each $\Primes^{\dag}$-integer $N>0$, and 
$$\Sha (A,\Cal C)\defeq \Ker \lgroup H^1(G_K,A)\to \prod_{c\in \Cal C^1} H^1(G_{K_c},A_c)\rgroup.$$
Thus, we have 
$$\Sel_N (A,\Cal C)\subset \Sel_N(A)=\Sel_N (A,C)$$ 
and 
$$\Sha (A,\Cal C)\subset \Sha(A)=\Sha (A,C).$$
In [Lang-Tate], Theorem 3 and Theorem 5, it is shown that 
$\Sel_N (A,\Cal C)$ and $\Sha(A,\Cal C)[N]$ are finite. Thus, 
Proposition 2.10 and Proposition 3.9 (i) 
can be regarded as an improvement of these results, respectively. 
\enddefinition

\subhead
\S4. The isotrivial case 
\endsubhead
We use the same notations as in $\S1$, $\S2$ and $\S3$. 
In particular, 
$k$ is a field which is {\bf finitely generated} over the prime field of 
characteristic $p\geq 0$ and infinite; 
$C\to \Spec k$ is a smooth, separated and geometrically connected 
algebraic curve over $k$ with 
function field $K=k(C)$; and 
$\Cal A\to C$ is an abelian scheme over $C$ with 
generic fibre $A= \Cal A_K = \Cal A \times _C \Spec K$. 


\proclaim {Theorem 4.1} 
Assume that 
$A$ is {\bf essentially isotrivial}, i.e., 
$A_{\overline K}$ is isogenous to an abelian variety over $\overline K$ 
that descends to an abelian variety over $\overline k$. 
Then $\Sha(A)^{(\dag)}$ is {\bf finite}.
\endproclaim

\demo{Proof} 
By Proposition 3.12 and Proposition 3.11, we may assume that $C(k)\nemp$, that 
$C$ admits a (unique) smooth compactification $C^{\cpt}$, 
and that $A$ is constant, i.e., 
$A$ descends to an abelian variety $\widetilde A$ 
over $k$ ($\widetilde A\times _{k}K=A$). 
In order to prove that
$\Sha(A)^{(\dag)}$ is finite, it suffices to show that $\Sha(A)[l^{\infty}]$ is finite 
for all $l\in\Primes^{\dag}$ and that 
$\Sha (A)[l^{\infty}]=0$ for all but finitely many $l\in \Primes ^{\dag}$. 
Further, as $\Sha(A)[l^n]$ is finite for all $n\geq 0$ by Proposition 3.9, the condition 
that $\Sha(A)[l^{\infty}]$ is finite is equivalent to: $T_l\Sha(A)=0$. 

Let $l\in \Primes ^{\dag}$. 
We view $(T_l\widetilde A\isom)T_lA$, which is fixed by $\pi_1(C_{\overline k})$, as a $G_k$-module, and we identify $H^1(G_k,T_lA)$ with $H^1(G_k,T_l\widetilde A)$. 
We have a natural inflation-restriction exact sequence 
$$0\to H^1(G_k,T_l\widetilde A)@>\inf>>  H^1(\pi_1(C),T_lA) @>\res>> \Hom (\pi_1(C_{\overline k}),T_lA)^{G_k}.$$
First, observe that (cf. Definition 2.3 (ii) for the definition of $\Sel^{\{l\}}(A)$)
$$\Sel^{\{l\}}(A)\cap H^1(G_k,T_l\widetilde A)=
A(K)^{\wedge,l}\cap H^1(G_k,T_l\widetilde A)=\widetilde A(k)^{\wedge,l}$$ 
holds in $H^1(\pi_1(C),T_lA)$. Indeed, the inclusions 
$$\Sel^{\{l\}}(A)\cap H^1(G_k,T_l\widetilde A)\supset
A(K)^{\wedge,l}\cap H^1(G_k,T_l\widetilde A)\supset \widetilde A(k)^{\wedge,l}$$ 
are clear. To prove
$\Sel^{\{l\}}(A)\cap H^1(G_k,T_l\widetilde A)\subset \widetilde A(k)^{\wedge,l}$, 
fix $c\in C(k)\nemp$. Then the composite of 
the inflation map 
$H^1(G_k,T_l\widetilde A) \to H^1(\pi_1(C),T_lA)$ 
and the restriction map 
$H^1(\pi_1(C),T_lA)\to H^1(G_{k(c)},T_l\Cal A_c)=H^1(G_k,T_l\widetilde A)$ 
at $c$ is the identity. As the image of $\Sel^{\{l\}}(A)$ 
under the restriction map 
$H^1(\pi_1(C),T_lA)\to H^1(G_{k(c)},T_l\Cal A_c)=H^1(G_k,T_l\widetilde A)$ 
at $c$ is included in $\widetilde A(k)^{\wedge,l}$ by definition, 
we obtain the desired inclusion. 

Now, let $\varphi :A(K)^{\wedge,l} \to \Hom (\pi_1(C_{\overline k}),T_lA)^{G_K}$ be the composite of the natural maps
$A(K)^{\wedge,l} \hookrightarrow H^1(\pi_1(C),T_lA)\to \Hom (\pi_1(C_{\overline k}),T_lA)^{G_K}$. 
Thus, we have a natural map 
$$T_l\Sha(A)={\Sel^{\{l\}}(A)}/
A(K)^{\wedge,l}\to \Hom (\pi_1(C_{\overline k}),T_lA)^{G_k}/\varphi (A(K)^{\wedge,l}),$$
which is injective as 
$$\Sel^{\{l\}}(A)\cap H^1(G_k,T_l\widetilde A)=\widetilde A(k)^{\wedge,l}\subset A(K)^{\wedge,l}.$$ 
To prove that 
$T_l\Sha(A)=0$, 
it suffices to show that 
$\Hom (\pi_1(C_{\overline k}),T_lA)^{G_k}/\varphi (A(K)^{\wedge,l})=0$. 
The Tate conjecture for abelian varieties holds over finitely generated fields by Tate, Zarhin, Mori 
in positive characteristic and by Faltings in characteristic $0$ 
(cf. [Tate], [Zarhin1], [Moret-Bailly] and [Faltings1]). As a consequence, we have 
a natural isomorphism 
$$\Hom_k(J,\widetilde A)\otimes _{\Bbb Z}\Bbb Z_l
\isom 
\Hom(T_lJ,T_l\widetilde A)^{G_k} 
=
\Hom (\pi_1(C^{\cpt}_{\overline k})^{\ab,l},T_l\widetilde A)^{G_k},$$
where $J$ denotes the jacobian variety of $C^{\cpt}$. 
We also have natural isomorphisms 
$$\Hom (\pi_1(C^{\cpt}_{\overline k})^{\ab,l},T_l\widetilde A)^{G_k}
\isom 
\Hom (\pi_1(C_{\overline k})^{\ab,l},T_l\widetilde A)^{G_k}
\isom 
\Hom (\pi_1(C_{\overline k}),T_l\widetilde A)^{G_k}, 
$$ 
where the first isomorphism follows from a standard weight argument. More precisely, we can assume without loss of generality 
(after possibly replacing $k$ by a finite extension)
that $C^{\cpt}\setminus C=\{c_0,c_1,\ldots,c_n\}\subset C^{\cpt}(k)$. Thus, 
$I_C\defeq 
I_C^{(\ab,l)}\defeq 
\Ker  \left(\pi_1(C_{\overline k})^{\ab,l}\twoheadrightarrow  \pi_1(C^{\cpt}_{\overline k})^{\ab,l}\right)
\isom \Coker \left (\Bbb Z_l(1)@>\diag>> \oplus_{i=0}^n\Bbb Z_l(1)\right)\isom
\oplus_{i=1}^n\Bbb Z_l(1)$ as $G_k$-module. Now the $G_{k}$-representation
$I_C\otimes \Bbb Q_l$ (resp. $T_l\widetilde A\otimes \Bbb Q_l$) is pure of weight $-2$ (resp. pure of weight $-1$) (cf. [Jannsen], 2), and
$\Hom(I_C, T_l\widetilde {\Cal A})^{G_{k}}=0$ follows (cf. loc. cit. Fact 2). Hence $\Hom (\pi_1(C^{\cpt}_{\overline k})^{\ab,l},T_l\widetilde A)^{G_k}
\isom 
\Hom (\pi_1(C_{\overline k})^{\ab,l},T_l\widetilde A)^{G_k}$. Further, 
as $C^{\cpt}(k)\nemp$, the natural map 
$$A(K)=\Mor_k(C^{\cpt},\widetilde A)\to \Hom_k(J,\widetilde A)$$ 
induced by the Albanese property of $J$ is surjective. 
Thus, the above map $\varphi :A(K)^{\wedge,l} \to \Hom (\pi_1(C_{\overline k}),T_l\widetilde A)^{G_k}$ 
is surjective and  
$\Hom (\pi_1(C_{\overline k}),T_l\widetilde A)^{G_k}/\varphi (A(K)^{\wedge,l})$ is trivial, as desired. 

Next, we prove that $\Sha (A)[l^{\infty}]=0$, or, equivalently, 
$\Sha (A)[l]=0$, for all 
but finitely many $l\in \Primes ^{\dag}$. Indeed, this follows from a similar argument as above
using the following truncated version of the Tate conjecture 
$$\Hom_k(J,\widetilde A)\otimes _{\Bbb Z}\Bbb Z/l\Bbb Z
\isom 
\Hom(J[l],\widetilde A[l])^{G_k}
,$$
which holds for all but finitely many $l\in\Primes^{\dag}$ 
(cf. [Zarhin2], [Zarhin3] in positive characteristic and [Faltings2], VI, \S3 
in characteristic $0$). 
\qed
\enddemo

\proclaim {Theorem 4.2} 
Assume that 
$A$ is {\bf essentially isotrivial}. 
Then the assertion of Conjecture 3.8 holds (resp. holds up to $p$-torsion) 
if $p=0$ (resp. $p>0$). More precisely, we have $\Sh^{\Sigma}(A)^{(\dag)}=0$.
\endproclaim

\demo{Proof} 
This follows immediately from Theorem 4.1 and Proposition 3.7. 
\qed
\enddemo

$$\text{References.}$$

\noindent
[Faltings1] Faltings, G., Endlichkeitss\"atze  f\"ur Abelsche Variet\"aten \"uber Zahlk\"orpern, Inventiones Math., 73 (1983), 349--366.

\noindent
[Faltings2] Faltings, G., Complements to Mordell, in Rational Points, Seminar Bonn/Wuppertal 1983/84, 
Faltings, G., 
W\"ustholz, G. 
et al., Aspects of Mathematics, E6, Third enlarged edition, 
Vieweg, 
1992.

\noindent
[Grothendieck] Grothendieck, A., Rev\^etements \'etales et groupe fondamental, Lecture 
Notes in Math., 224, Springer, 
1971.

\noindent
[Jannsen] Jannsen, U., Weights in arithmetic geometry, Japanese Journal of Mathematics, 
5 (2010), no. 1, 
73--102.

\noindent
[Katz-Lang] Katz, N. M., Lang, S, Finiteness theorems in geometric classfield theory, 
Enseign. Math. (2),  27  (1981), no. 3-4, 285--319.

\noindent
[Lang-N\'eron] Lang, S., N\'eron, A, Rational points of abelian varieties over function fields, Amer. J. Math., 81 (1959), 95--118.

\noindent
[Lang-Tate] Lang, S., Tate, J., Principal homogeneous spaces over abelian varieties, Amer. J. Math., 80 (1958), 659--684. 

\noindent
[Moret-Bailly] Moret-Bailly, S., Pinceaux de Vari\'et\'es Ab\'eliennes, Ast\'erisque, 129, Soc. Math. France, 1985. 

\noindent
[Neukirch-Schmidt-Wingberg] Neukirch, J., Schmidt, A., Wingberg, K., Cohomology of Number Fields, Grundlehren der mathematischen Wissenschaften, 323, Springer, 2000.

\noindent
[Ribes-Zalesskii] Ribes, L., Zalesskii, P., Profinite Groups, Ergebnisse der Mathematik und ihrer Grenzgebiete, 
Folge 3, 
40, Springer, 2000.

\noindent
[Sa\"\i di1] Sa\"\i di, M., The cuspidalisation of sections of arithmetic fundamental groups, 
Advances in Mathematics, 230 (2012), 1931--1954.

\noindent
[Sa\"\i di2] Sa\"\i di, M., On the section conjecture over function fields and finitely generated fields, Pub. Res. Inst. Math. Sci., 
52 (2016), no. 3, 335-357.

\noindent
[Serre1] Serre, J.-P., Galois cohomology, Springer-Verlag Berlin Heidelberg, 1997.

\noindent
[Serre2] Serre, J.-P., Lectures on the Mordell-Weil Theorem, Translated and edited by 
Brown, M., from notes by 
Waldschmidt, M., 
Second edition, 
Aspects of Mathematics, 
Vieweg, 
1990.

\noindent
[Serre-Tate] Serre, J.-P., Tate, J., Good reduction of abelian varieties, Annals of Mathematics,
Second series, 88 (1968), no. 3, 492--517.

\noindent
[Tate] Tate, J., Endomorphisms of abelian varieties over finite fields, Inventiones Math., 2 (1966), 134--144.


\noindent
[Zarhin1] Zarkhin, Yu. G., Abelian varieties in characteristic $p$, Math. Notes, 19 (1976), no. 3, 240--244. 

\noindent
[Zarhin2] Zarkhin, Yu. G., Endomorphisms of Abelian varieties and points of finite order in characteristic $p$, 
Math. Notes, 21 (1977), no. 6, 415--419. 

\noindent
[Zarhin3] Zarhin, Yu. G., Abelian varieties over fields of finite characteristic, Cent. Eur. J. Math., 12 (2014), no. 5, 
659--674. 

\bigskip
\noindent
Mohamed Sa\"\i di

\noindent
College of Engineering, Mathematics, and Physical Sciences

\noindent
University of Exeter

\noindent
Harrison Building

\noindent
North Park Road

\noindent
EXETER EX4 4QF 

\noindent
United Kingdom

\noindent
M.Saidi\@exeter.ac.uk

\bigskip
\noindent
Akio Tamagawa

\noindent
Research Institute for Mathematical Sciences

\noindent
Kyoto University

\noindent
KYOTO 606-8502

\noindent
Japan

\noindent
tamagawa\@kurims.kyoto-u.ac.jp
\enddocument